\documentclass{amsart}
\usepackage[foot]{amsaddr}

\usepackage{amssymb}
\usepackage{cases}
\usepackage{tikz}
\usepackage{anyfontsize}
\usepackage{enumitem}
\setlist[itemize]{leftmargin=2em}

\newcommand\bR{{\mathbb R}}

\newcommand\bB{{\mathbb B}}
\newcommand\gph{{\mathrm{gph}\,}}
\newcommand\ior{{\mathrm{int}\,}}
\newcommand\ri{{\mathrm{ri}\,}}
\newcommand\rb{{\mathrm{rb}\,}}

\newcommand\sgn{{\mathrm{sgn}\,}}
\newcommand\dist{{\mathrm{dist}}}

\newcommand\ind{{\mathrm{ind}}}

\numberwithin{equation}{section}

\newtheorem{theorem}{Theorem}[section]
\newtheorem{corollary}[theorem]{Corollary}
\newtheorem{lemma}[theorem]{Lemma}
\newtheorem{proposition}[theorem]{Proposition}
\newtheorem{remark}[theorem]{Remark}

\usepackage[sort&compress]{natbib}
 \bibpunct[, ]{[}{]}{,}{n}{}{,}%

\usepackage[colorlinks=true,breaklinks=true,bookmarks=true,urlcolor=blue,citecolor=blue,linkcolor=blue,bookmarksopen=false,draft=false]{hyperref}

\usepackage{xparse}
\makeatletter
\RenewDocumentCommand{\title}{om}{%
   \IfNoValueTF{#1}
     {\gdef\shorttitle{The Aubin Property for Generalized Equations}}
     {\gdef\shorttitle{#1}}%
   \gdef\@title{#2}%
}
\makeatother

\title{The Aubin Property for Generalized Equations over $C^2$-cone Reducible Sets}
\author{Jiaming Ma}
\email{jiaming.ma@connect.polyu.hk}
\author{Defeng Sun}
\email{defeng.sun@polyu.edu.hk}
\address{Department of Applied Mathematics, The Hong Kong Polytechnic University}
\date{August 10, 2026}

\thanks{This research was supported in part by the Hong Kong Research Grants Council under GRF project 15309625 and the RGC Senior Research Fellow scheme SRFS2223-5S02.}

\begin{document}

\begin{abstract}
    This paper establishes the equivalence between the Aubin property and strong regularity for generalized equations over $C^2$-cone reducible sets. This result resolves a long-standing question in variational analysis and extends the classical equivalence theorem for polyhedral sets to a broad class of non-polyhedral sets. Our proof strategy departs from traditional variational techniques, integrating insights from convex geometry with powerful tools from algebraic topology. At the heart of our analysis is a novel index theorem for a class of functions involving metric projections onto arbitrary closed convex sets. Its proof exploits the geometry of normal cones and topological degree theory. We then use a lift of the diffeomorphism provided by the $C^2$-cone reduction to reduce the original generalized equations to the setting of the index theorem. The homological inverse mapping theorem then yields the desired strong regularity. {We show that the $C^2$-cone reducibility assumption cannot be removed in general by constructing a three-dimensional semialgebraic counterexample.} This result unifies and extends existing stability results for conventional nonlinear programming, nonlinear second-order cone programming, and nonlinear semidefinite programming under a single general framework, and removes the local optimality assumption imposed in the existing non-polyhedral results.

    \smallskip \noindent \textbf{Keywords.} Aubin property, Strong regularity, $C^2$-cone reducible set, Generalized equation, Convex geometry, Degree theory
\end{abstract}

\maketitle

\section{Introduction}

In this paper, we establish that the Aubin property and strong regularity are equivalent for the solution mappings of canonically perturbed generalized equations over $C^2$-cone reducible sets.
Specifically, we consider the following two types of generalized equations:
\begin{equation}\label{ge}
    y\in \varphi(x)+N_S(x)
    \quad\mathrm{and}\quad
    y\in \varphi(x)+N^{-1}_S(x), \quad x \in \bR^n, 
\end{equation}
where $\varphi:\bR^n\to\bR^n$ is a continuously differentiable function, $S\subset\bR^n$ is a nonempty closed convex set, $y\in\bR^n$ is the parameter vector and $N_S:\bR^n\rightrightarrows\bR^n$ is the normal cone map of $S$. Many problems in optimization and variational analysis can be written in the form of (\ref{ge}), such as the Karush-Kuhn-Tucker (KKT) system and the variational inequality. For convenience, let $\Phi:\bR^n\rightrightarrows\bR^n$ denote the right-hand side of the generalized equation in (\ref{ge}), i.e.,
\begin{equation}\label{Phi}
    \Phi(x)=\varphi(x)+N_S(x)
    \quad\mathrm{or}\quad
    \Phi(x)=\varphi(x)+N^{-1}_S(x), \quad x \in \bR^n.
\end{equation}

In a seminal paper \cite{dontchev1996characterizations}, Dontchev and Rockafellar established the equivalence between the Aubin property and strong regularity for $\Phi^{-1}$ when $S$ is polyhedral. Their argument, developed for the first form of $\Phi$, relies essentially on the polyhedral structure of $S$ and the associated piecewise affine normal map introduced by Robinson \cite{robinson1992normal}. Although the second form can also be reduced to the polyhedral framework, this approach does not extend directly to non-polyhedral sets.

More recently, Chen et al. established the same equivalence for the KKT solution mappings of canonically perturbed nonlinear second-order cone programming \cite{chen2025aubin} and nonlinear semidefinite programming \cite{chen2025characterizations} at a locally optimal solution. These results concern important non-polyhedral cones, but their proofs make essential use of the specific geometry of the second-order cone and the positive semidefinite cone. Moreover, the local optimality assumption distinguishes these results from the classical polyhedral theorem.

The sets arising in all these results are $C^2$-cone reducible in the sense of \cite[Definition 3.135]{bonnans2000perturbation}. This naturally raises the question of whether the Aubin property and the strong regularity for $\Phi^{-1}$ remain equivalent for an arbitrary $C^2$-cone reducible set $S$. The main theorem of this paper gives an affirmative answer for both forms of $\Phi$ in \eqref{Phi}.

For a set-valued map $\Psi:\bR^n\rightrightarrows\bR^m$, the inverse of $\Psi$ is defined as $\Psi^{-1}(y):=\{x\in\bR^n\mid y\in\Psi(x)\}$. The graph of $\Psi$ is $\gph\Psi:=\{(x,y)\in\bR^n\times\bR^m\mid y\in\Psi(x)\}$. We are concerned with the following two localized Lipschitzian properties (the Aubin property and strong regularity) for $\Psi$ around $(x_0,y_0)\in\gph \Psi$:
\begin{enumerate}[label=\textbf{(L\arabic*)},leftmargin=3em]
    \item\label{L1} $\Psi$ has the Aubin property \cite{aubin1984lipschitz} around $(x_0,y_0)$. That is, there exist open neighborhoods $U$ of $x_0$ and $V$ of $y_0$ and a positive constant $\lambda>0$ such that
    \begin{equation*}
        \Psi(x)\cap V\subset\Psi(x^\prime)+\lambda\|x-x^\prime\|\bB,\quad\forall\, x, x^\prime\in U,
    \end{equation*}
    where $\bB$ is the closed unit ball in $\bR^m$.
    \item\label{L2} $\Psi$ is locally single-valued and Lipschitz continuous around $(x_0,y_0)$. That is, there exist open neighborhoods $U$ of $x_0$ and $V$ of $y_0$ such that the map $x\mapsto\Psi(x)\cap V$ is single-valued and Lipschitz continuous on $U$.
\end{enumerate}

Let $(x_0,y_0)\in\gph\Phi$. In a landmark paper \cite{robinson1980strongly}, Robinson found that generalized equations $y\in\Phi(x)$ are closely related to the corresponding linearized generalized equations $y\in\widetilde{\Phi}(x)$, where
\begin{equation*}
    \widetilde\Phi(x)=\varphi(x_0)+\varphi^\prime(x_0)(x-x_0)+N_S(x)
    \quad\mathrm{or}\quad
    \widetilde\Phi(x)=\varphi(x_0)+\varphi^\prime(x_0)(x-x_0)+N^{-1}_S(x).
\end{equation*}
The property \ref{L2} for $\widetilde\Phi^{-1}$ around $(y_0,x_0)$ is called the strong regularity for $\Phi^{-1}$ around $(y_0,x_0)$; see \cite{robinson1980strongly} and \cite[Definition 5.12]{bonnans2000perturbation}. Moreover, according to \cite[Theorem 5.13]{bonnans2000perturbation}, the property \ref{L2} for $\Phi^{-1}$ around $(y_0,x_0)$ is equivalent to the property \ref{L2} for $\widetilde{\Phi}^{-1}$ around the same point. Consequently, in the present setting, strong regularity of $\Phi^{-1}$ around $(y_0,x_0)$ is equivalent to the local single-valuedness and Lipschitz continuity of $\Phi^{-1}$ there.

Clearly, for any set-valued map, \ref{L2} implies \ref{L1}, but the converse does not necessarily hold in general. In this paper, we prove the equivalence between the properties \ref{L1} and \ref{L2} for $\Phi^{-1}$ around $(y_0,x_0)\in\gph\Phi^{-1}$ when $S$ is $C^2$-cone reducible. The $C^2$-cone reducibility assumption is sharp in general. We construct a three-dimensional closed convex semialgebraic set $S$, and not $C^2$-cone reducible at the reference point, together with a generalized equation of the form $\Phi=\nabla f+N_S$, such that $\Phi^{-1}$ has the Aubin property but is not strongly regular. In this example, the localized solution mapping is two-valued. Thus, the $C^2$-cone reducibility assumption cannot in general be removed from our main theorem.

As an application, we consider the canonically perturbed KKT system of a general constrained optimization problem. We prove that the Aubin property and strong regularity of its KKT solution mapping are equivalent at an arbitrary KKT point whenever the constraint set is $C^2$-cone reducible, or is a cone whose polar is $C^2$-cone reducible. In particular, this removes the local optimality assumption imposed in the existing nonlinear second-order cone and semidefinite programming results. The resulting framework also covers Cartesian products involving polyhedral cones, positive semidefinite cones, and broad classes of $p$-order cones.

From now on, we always assume that $S$ in \eqref{ge} is $C^2$-cone reducible. By the $C^2$-cone reduction and other operations, we find that the local behavior of $\Phi$ around $(x_0,y_0)$ is closely related to that of a function $N:\bR^n\to\bR^n$ around a point $z_0\in K$. $N$ is defined in the following form:
\begin{equation}\label{N}
N(z)=A(z-\Pi_K(z))+B\Pi_K(z),\quad z\in \bR^n,
\end{equation}
where $K\subset\bR^n$ is a nonempty closed convex set, $A$ and $B$ are $n$ by $n$ matrices, and $\Pi_K:\bR^n\to\bR^n$ is the metric projection onto $K$ given by
$$
\Pi_K(z):=\arg\min_{z^\prime\in K}\|z-z^\prime\|,\quad z\in\bR^n.
$$
Note that $N$ is similar to the normal map
$$
z\mapsto\varphi(\Pi_S(z))+(z-\Pi_S(z))
$$
induced by the generalized equation $y\in\varphi(x)+N_S(x)$ and proposed by Robinson \cite{robinson1992normal}. In fact, with appropriately chosen $K$, $A$, and $B$, the local properties of $\Phi$ can be related to those of $N$ around some point $z_0\in K$. The condition $z_0\in K$ is crucial for our topological approach, and the $C^2$-cone reduction enables us to obtain such a representation.

Since the set $K$ in \eqref{N} is an arbitrary closed convex set without any additional specific structure, we shall use a topological approach combined with convex geometry. Our proof relies on two key tools: the geometry of the normal cone map $N_K$, particularly the dimensions of its values, and topological degree theory \cite{lloyd1978degree,dinca2021brouwer}. Degree theory is a powerful tool for analyzing the existence, multiplicity, and qualitative properties of solutions to equations involving continuous functions, even in the absence of differentiability. It is therefore suitable
for the study of the continuous function $N$ and the normal map.

The general idea of our proof goes as follows. Given $\Phi$ and $(x_0,y_0)\in\gph\Phi$, suppose that $\Phi^{-1}$ has the Aubin property around $(y_0,x_0)$. There is a continuous function $\widetilde{N}:\bR^n\to\bR^n$, similar to a normal map and induced by $\Phi$, such that $\widetilde{N}^{-1}$ also has the Aubin property around $(\widetilde{N}(\tilde{z}_0),\tilde{z}_0)$, where $\tilde{z}_0\in\bR^n$ is determined by $\Phi$ and $(x_0,y_0)$. By the $C^2$-cone reduction and suitable local transformations, the behavior of $\widetilde{N}$ around $\tilde{z}_0$ can be related to that of a function $N$ of the form \eqref{N} around a point $z_0\in K$.

The most critical step is to use the geometric structure of $N_K$ around $z_0$ to prove that the topological index of $N$ at $z_0$ is $\pm1$. The required openness and discreteness conditions are ultimately obtained from the Aubin property of $\widetilde{N}^{-1}$ and the directional differentiability of $\widetilde{N}$. The index information can then be transferred back to $\widetilde{N}$, showing that its index at $\tilde{z}_0$ is also $\pm1$. According to the homological version of the inverse mapping theorem \cite{barreto2016inverse}, $\widetilde{N}$ is a local homeomorphism at $\tilde{z}_0$, which yields the strong regularity of $\Phi^{-1}$ around $(y_0,x_0)$. The detailed proofs are presented in the following sections.

The remainder of this paper is organized as follows. Section \ref{sec:pre} introduces the notation and preliminary results concerning convex geometry, degree theory, strictly stationary perturbations, and $C^2$-cone reductions. In Section \ref{sec:index}, we establish the index theorem for the canonical projection-based map over an arbitrary closed convex set. Section \ref{sec:eqv} proves the equivalence between the Aubin property and strong regularity for the generalized equations over $C^2$-cone reducible set. {Section \ref{sec:counter} constructs a counterexample showing that the $C^2$-cone reducibility assumption cannot be removed in general. Section \ref{sec:app} applies the main theorem to canonically perturbed KKT systems.} Finally, Section \ref{sec:con} concludes the paper.

\section{Notation and preliminaries}\label{sec:pre}

We use the standard inner product $\langle\cdot,\cdot\rangle$ and the Euclidean norm $\|\cdot\|$ in $\bR^n$.
Let $\Omega\subset\bR^n$. We use $\ior\Omega$ and $\overline{\Omega}$ to denote the interior and closure of $\Omega$, respectively. The boundary of $\Omega$ is denoted by $\partial\Omega:=\overline{\Omega}\backslash\ior\Omega$. 
A point $x\in\Omega$ is an isolated point of $\Omega$ if there exists a neighborhood $U$ of $x$ such that $U\cap\Omega=\{x\}$.
The distance from $x\in\bR^n$ to $\Omega$ is defined by $\dist(x,\Omega):=\inf\{\|x-x^\prime\|\mid x^\prime\in\Omega\}$.
Let $L\subset\bR^n$ be a linear subspace. The orthogonal complement of $L$ is denoted by $L^\perp$. The orthogonal projection matrix $P$ onto $L$ is an $n$ by $n$ matrix such that $Px=x$ for all $x\in L$ and $Px=0$ for all $x\in L^\perp$.

For a nonempty closed convex cone $C\subset\bR^n$, the polar cone of $C$ is given by $C^\circ:=\{u\in\bR^n\mid\langle x,u\rangle\leq 0,\:\forall x\in C\}$. Let $K\subset\bR^n$ be a nonempty closed convex set. The normal cone map $N_K:\bR^n\rightrightarrows\bR^n$ is defined by
\begin{equation*}
    N_K(x):=
    \begin{cases}
        \{u\in\bR^n\mid\langle u,x^\prime-x\rangle\leq0,\:\forall x^\prime\in K\} & \text{if } x \in K, \\
        \emptyset & \text{if } x \notin K,
    \end{cases}
\end{equation*}
and the tangent cone map $T_K:\bR^n\rightrightarrows\bR^n$ of $K$ is given by
\begin{equation*}
    T_K(x):=
    \begin{cases}
        N_K^\circ(x) & \text{if } x \in K, \\
        \emptyset & \text{if } x \notin K.
    \end{cases}
\end{equation*}
Note that for any $x\in K$, $u\in N_K(x)$ if and only if $x=\Pi_K(x+u)$. The dimension of $K$, i.e., the dimension of the affine hull of $K$, is denoted by $\dim K$. The relative interior and relative boundary of $K$ are denoted by $\ri K$ and $\rb K$, respectively. A face of $K$ is a convex subset $F\subset K$ such that $\lambda x+(1-\lambda)x^\prime\in F$ with $x,x^\prime\in K$ and $0<\lambda<1$ implies $x,x^\prime\in F$. A supporting half-space to $K$ is a closed half-space containing $K$ and having a point of $K$ in its boundary, which can be represented in the form $\{x\in\bR^n\mid\langle x,u\rangle\leq\alpha\}$ where $u\neq 0$ is called an (outer) normal vector of $K$. A supporting hyperplane to $K$ is the boundary of a supporting half-space to $K$. For each supporting hyperplane $H$ to $K$, the set $K\cap H$ is a face of $K$ called an exposed face of $K$.

A continuously differentiable function is called the $C^1$ function, and a twice continuously differentiable function is called the $C^2$ function. If $f:\bR^n\to\bR^m$ is differentiable at $x_0\in\bR^n$, the Jacobian matrix of $f$ at $x_0$ is denoted by $f'(x_0)$, and we write $\nabla f(x_0):=(f'(x_0))^T$. Let $U,V\subset\mathbb R^n$ be open sets. A $C^k$ homeomorphism $h:U\to V$ is called a $C^k$-diffeomorphism if $\det h'(x)\neq 0$ for every $x\in U$.

Let $f:\bR^n\to\bR^m$ be a function and let $x_0\in\bR^n$. We say that $f$ is a local homeomorphism at $x_0$ if there exists an open neighborhood $U$ of $x_0$ such that $f(U)$ is open and $f|_U:U\to f(U)$ is a homeomorphism. If $m=n$, we say that $f$ is a local $C^k$-diffeomorphism at $x_0$ if $f|_U\to f(U)$ is a $C^k$-diffeomorphism for some open neighborhood $U$ of $x_0$.

We say that $f$ is open around $x_0$, if there exists an open neighborhood $U$ of $x_0$ such that $f|_U$ (the restriction of $f$ to $U$) is an open map, that is, for any open set $V\subset U$, $f(V)$ is open in $\bR^m$. The function $f$ is said to be discrete at $x_0$ if $x_0$ is an isolated point of $f^{-1}(f(x_0))$, i.e., there exists an open neighborhood $W$ of $x_0$ such that $f^{-1}(f(x_0))\cap W=\{x_0\}$. If there exists an open neighborhood $U$ of $x_0$ such that $f$ is discrete at each $x\in U$, then we say that $f$ is discrete around $x_0$.

Given a set-valued map $\Psi:\bR^n\rightrightarrows\bR^m$ and $x_0\in\bR^n$. $\Psi$ is outer semicontinuous at $x_0$ if
\begin{equation*}
    \limsup_{x\to x_0} \Psi(x)\subset\Psi(x_0)
\end{equation*}
and $\Psi$ is inner semicontinuous at $x_0$ if
\begin{equation*}
    \liminf_{x\to x_0} \Psi(x)\supset\Psi(x_0),
\end{equation*}
where $``\limsup"$ and $``\liminf"$ are the outer limit and the inner limit in Painlevé-Kuratowski convergence for subsets, respectively.
$\Psi$ is called outer (inner) semicontinuous if $\Psi$ is outer (inner) semicontinuous at each $x\in\bR^n$. According to \cite[Corollary 6.29]{rockafellar1998variational}, for a nonempty closed convex set $K\subset\bR^n$, the normal cone map $N_K:\bR^n\rightrightarrows\bR^n$ is outer semicontinuous, and the tangent cone map $T_K:\bR^n\rightrightarrows\bR^n$ is inner semicontinuous. We say $\Psi$ is inner semicontinuous around $(x_0,y_0)\in\gph\Psi$ if there exist open neighborhoods $U$ of $x_0$ and $V$ of $y_0$ such that the map $x\mapsto\Psi(x)\cap V$ is inner semicontinuous at each $x\in U$. We have the following lemma.
\begin{lemma}\label{open}
    Let $f:\bR^n\to\bR^m$ be a continuous function and $x_0\in\bR^n$. If $f^{-1}$ is inner semicontinuous around $(f(x_0),x_0)$, then $f$ is open around $x_0$. 
\end{lemma}
\begin{proof}
    Assume that $f^{-1}$ is inner semicontinuous around $(f(x_0),x_0)$. Then there exist open neighborhoods $U$ of $x_0$ and $V$ of $f(x_0)$ such that the following set-valued map $\Psi:\bR^m\rightrightarrows\bR^n$ is inner semicontinuous:
    $$
    \Psi(y):=
    \begin{cases}
        \{x\in U\mid f(x)=y\} & \text{if } y \in V, \\
        \emptyset & \text{if } y \notin V.
    \end{cases}
    $$
    So $\Psi^{-1}(W)$ is open for every open set $W\subset \bR^n$ \cite[Theorem 5.7, Part (c)]{rockafellar1998variational}. Let $W\subset U\cap f^{-1}(V)$ be an open set. Then the set $\Psi^{-1}(W)=\{y\in\bR^n\mid\Psi(y)\cap W\neq\emptyset\}=\{y\in V\mid \{x\in W\mid f(x)=y\}\neq\emptyset\}=f(W)$ is open. Therefore, $f$ is open around $x_0$, since $U\cap f^{-1}(V)$ is an open set containing $x_0$.
\end{proof}

We know, by definition, that the Aubin property of $\Psi$ around $(x_0,y_0)\in\gph\Psi$ implies that $\Psi$ is inner semicontinuous around $(x_0,y_0)$. So, by Lemma \ref{open}, for a continuous function $f:\bR^n\to\bR^m$ and $x_0\in\bR^n$, if $f^{-1}$ has the Aubin property around $(f(x_0),x_0)$, then $f$ is open around $x_0$.

\subsection{Classification of normal vectors}\label{sec:nv}

Consider the function $N$ defined in (\ref{N}). When $B=I$, we have
\begin{equation*}
    N(x+u)=A(x+u-\Pi_K(x+u)) + \Pi_K(x+u) = Au +x , \quad \forall\, (x,u)\in \gph N_K.
\end{equation*}
So in a sense, the function $N$ acts as a ``linear transform'' on the normal vectors of $K$. Therefore, it is important to understand the structure of the normal vectors of $K$ around a point $x\in K$.

We need the classification of normal vectors of a compact convex set $K\subset\bR^n$ described in \cite[Section 2.2]{schneider2013convex}, which was first studied by Bonnesen and Fenchel \cite{BonnesenFenchel1934}. Let $u\in\bR^n\backslash\{0\}$. Then $N_K^{-1}(u):=\{x\in K\mid u\in N_K(x)\}$ is the exposed face of $K$ with outer normal vector $u$. For each $x\in \ri N_K^{-1}(u)$, the normal cone of $K$ at $x$ is the same; we denote this common normal cone by $N_K(\ri N_K^{-1}(u))$. The touching cone of $K$ at $u$, denoted by $T(K,u)$, is the smallest face of $N_K(\ri N_K^{-1}(u))$ that contains $u$. Moreover, for each $x\in N_K^{-1}(u)$, $T(K,u)$ is the smallest face of $N_K(x)$ that contains $u$; see \cite[Note 6 for Section 2.2]{schneider2013convex}. The vector $u$ is called an $r$-extreme normal vector of $K$ if $\dim T(K,u)\leq r+1$. And $u$ is called an $r$-exposed normal vector of $K$ if $\dim N_K(\ri N_K^{-1}(u))\leq r+1$. {Since $T(K,u)$ is a face of $N_K(\ri N_K^{-1}(u))$, every $r$-exposed normal vector is $r$-extreme. The following lemma, which is a dual version of the classical result of Asplund \cite{asplund1963k}, states that the converse holds after taking limits.}

\begin{lemma}[{\cite[Theorem 2.2.9]{schneider2013convex}}]\label{extexpthm}
    Let $K\subset\bR^n$ be a nonempty compact convex set and $r\in\{0,1,...,n-1\}$. Then each $r$-extreme normal vector of $K$ is a limit of $r$-exposed normal vectors of $K$.
\end{lemma}

\subsection{Degree theory and inverse mapping theorem}
For any bounded open set $D\subset\bR^n$, continuous function $f:\overline{D}\to\bR^n$, and $y\in\bR^n\backslash f(\partial D)$, the degree of $f$ on $D$ at $y$ is an integer denoted by $\deg(f,D,y)$. {For later use, we recall the following standard properties of the degree, collected from \cite[Section 1.2]{dinca2021brouwer}.}
\begin{proposition}
    Let $D\subset\bR^n$ be a bounded open set and $f:\overline{D}\to\bR^n$ be a continuous function.
    \begin{itemize}
        \item Let $y\in\bR^n\backslash f(\partial D)$. If $f$ is differentiable on an open set containing $f^{-1}(y)$ and the Jacobian matrix $f^\prime(x)$ is nonsingular for all $x\in f^{-1}(y) $, then $\deg(f,D,y)=\sum_{x\in f^{-1}(y)}\sgn \det f'(x)$.
        \item (Local constancy) The degree $\deg(f,D,\cdot)$ is constant on every connected component of $\bR^n\backslash f(\partial D)$.
        \item (Homotopy invariance) If $H:[0,1]\times\overline{D}\to\bR^n$ is continuous, $\gamma:[0,1]\to\bR^n$ is continuous, and $\gamma(t)\notin H(t,\partial D)$ for all $t\in[0,1]$, then $\deg(H(t,\cdot),D,\gamma(t))$ does not depend on $t$.
    \end{itemize}
\end{proposition}

Let $f:\bR^n\to\bR^n$ be a continuous function. Suppose that $x\in\bR^n$ is an isolated point of $f^{-1}(f(x))$. Take any bounded open set $D$ such that $f^{-1}(f(x))\cap D=\{x\}$. Then the degree $\deg(f|_{\overline{D}},D,f(x))$ is independent of the choice of $D$. This integer is called the index, or local degree, of $f$ at $x$, and is denoted by $\ind(f,x)$. {By the extended homotopy invariance property \cite[Theorem 2.1.1]{dinca2021brouwer}, for any $x',y'\in\bR^n$ we have $\ind(f,x)=\ind(f(\cdot+x^\prime)+y^\prime,x-x^\prime)$.}

We shall use the following homological version of the inverse mapping theorem \cite{barreto2016inverse}, which is equivalent to the homological version of the implicit function theorem \cite{biasi2008implicit}. This result is a consequence of the classical theorem for discrete open maps, proved by {\v{C}}ernavski{\u\i} \cite{chernavskii1964finite,chernavskii1965addendum} and later by V{\"a}is{\"a}l{\"a} \cite{vaisala1967discrete}.

\begin{lemma}[{\cite[Theorem 1.2]{barreto2016inverse}}]\label{localhom}
    Let $f:\bR^n\to\bR^n$ be a continuous function. If $f$ is open and discrete around $x_0\in\bR^n$ and $|\ind(f,x_0)|=1$, then $f$ is a local homeomorphism at $x_0$.
\end{lemma}

\subsection{Strictly stationary function}

A set-valued map $\Psi:\bR^n\rightrightarrows\bR^m$ is called locally closed-valued around $(x_0,y_0)\in\gph\Psi$, if there exist neighborhoods $U$ of $x_0$ and $V$ of $y_0$ such that $\Psi(x)\cap \overline{V}$ is closed for each $x\in U$. Consider the following property for $\Psi$ around $(x_0,y_0)\in\gph\Psi$:
\begin{enumerate}[label=\textbf{(L\arabic*)},leftmargin=3em]
    \setcounter{enumi}{2}
    \item\label{L3} $\Psi$ is locally closed-valued and has the Aubin property around $(x_0,y_0)$.
\end{enumerate}
Note that for a continuous function $f:\bR^n\to\bR^m$, properties \ref{L1} and \ref{L3} are the same for $f^{-1}$ since $f^{-1}(y)$ is closed for any $y\in\bR^m$.
A function $\psi:\bR^n\to\bR^m$ is called strictly stationary \cite{dontchev1994inverse} at $x_0\in\bR^n$ if, for any $\varepsilon>0$, there exists an open neighborhood $U$ of $x_0$, such that
\begin{equation*}
    \|\psi(x)-\psi(x^\prime)\|<\varepsilon\|x-x^\prime\|,\quad\forall\, x,x^\prime\in U.
\end{equation*}
\begin{lemma}[\cite{dontchev1994inverse}]\label{invmapthm}
    Let $\Psi:\bR^n\rightrightarrows\bR^m$ be a set-valued map and $(x_0,y_0)\in\gph\Psi$. Let $\psi:\bR^n\to\bR^m$ be a function that is strictly stationary at $x_0$. Then $\Psi^{-1}$ has the property \ref{L2} (respectively, \ref{L3}) around $(y_0,x_0)$ if and only if $(\Psi+\psi)^{-1}$ has the property \ref{L2} (respectively, \ref{L3}) around $(y_0+\psi(x_0),x_0)$.
\end{lemma}
Lemma \ref{invmapthm} shows that strictly stationary perturbations do not change the Aubin property of inverse mappings. So, it is important to find out what kind of functions are strictly stationary. Note that $\psi:\bR^n\to\bR^m$ is strictly stationary at $x_0$ if and only if each component $\psi_k:\bR^n\to\bR$ is strictly stationary at $x_0$. The following lemma records several elementary sufficient conditions for strictly stationary functions. The proof is direct and thus omitted.
\begin{lemma}\label{ss}
    Let $f:\bR^n\to\bR$, $g:\bR^n\to\bR$ and $x_0\in\bR^n$. Consider the following cases:
    \begin{itemize}
        \item $f,g$ are strictly stationary at $x_0$. Let $\psi(x)=f(x)+g(x)$.
        \item $f$ is strictly stationary at $y_0\in\bR^n$ and $h:\bR^n\to\bR^n$ is Lipschitz around $x_0$ such that $h(x_0)=y_0$. Let $\psi(x)=f(h(x))$.
        \item $f$ is continuously differentiable around $x_0$. Let $\psi(x)=f^\prime(x_0)x-f(x)$.
        \item $f$ is strictly stationary at $x_0$, $f(x_0)=0$ and $g$ is Lipschitz around $x_0$. Let $\psi(x)=f(x)g(x)$.
        \item $f,g$ are Lipschitz around $x_0$ and $f(x_0)=g(x_0)=0$. Let $\psi(x)=f(x)g(x)$.
    \end{itemize}
    Then $\psi$ is strictly stationary at $x_0$ in each case.
\end{lemma}

The following lemma allows us to keep the index under strictly stationary perturbations.
\begin{lemma}\label{keepindex}
    Let $f:\bR^n\to\bR^n$ be a continuous function such that $f$ is discrete at $x_0\in\bR^n$ and $f^{-1}$ has the Aubin property around $(f(x_0),x_0)$. Let $g:\bR^n\to\bR^n$ be a continuous function that is strictly stationary at $x_0$. Then $f+g$ is discrete at $x_0$ and $\ind(f,x_0)=\ind(f+g,x_0)$.
\end{lemma}
\begin{proof}
    Since $f^{-1}$ has the Aubin property around $(f(x_0),x_0)$ and $\gph f$ is closed, we have that $f$ is metrically regular at $(x_0,f(x_0))$ by \cite{borwein1988verifiable}, that is, there exist a constant $\kappa>0$ and bounded open neighborhoods $U$ of $x_0$ and $V$ of $f(x_0)$ such that for each $x\in U$ and $y\in V$, we have
    $$\dist(x,f^{-1}(y))\leq \kappa \dist(f(x),y).$$
    Since $f$ is discrete at $x_0$, by shrinking $U$ if necessary, we can assume that $f^{-1}(f(x_0))\cap U=\{x_0\}$ and for each $x\in \overline{U}$,
    $$\dist(x,f^{-1}(f(x_0)))=\|x-x_0\|\leq\kappa\|f(x)-f(x_0)\|.$$
    Since $g$ is strictly stationary at $x_0$, by shrinking $U$ if necessary, we can assume that for each $x\in \overline{U}$,
    $$\|g(x)-g(x_0)\|\leq (2\kappa)^{-1}\|x-x_0\|\leq\|f(x)-f(x_0)\|/2.$$
    Define the homotopy $H:[0,1]\times \overline{U}\to\bR^n$ by
    $$H(t,x):=f(x)+tg(x).$$
    Then $H(t,x_0)\notin H(t,\overline{U}\backslash\{x_0\})$ for all $t\in[0,1]$. By the homotopy invariance of the degree, we have
    $$\ind(f,x_0)=\deg(f|_{\overline{U}},U,f(x_0))=\deg((f+g)|_{\overline{U}},U,f(x_0)+g(x_0))=\ind(f+g,x_0), $$
    which completes the proof.
\end{proof}

Together, Lemmas \ref{invmapthm} and \ref{keepindex} allow us to add strictly stationary perturbations while preserving both the Aubin property of inverse and the index. This observation will be used in Section \ref{sec:eqv} to simplify the operators under consideration.

\subsection{\texorpdfstring{$C^2$}{C2}-cone reduction}

Let $S\subset\bR^n$ be a nonempty closed convex set. We say that $S$ is $C^2$-cone reducible at $x_0\in S$ \cite[Definition 3.135]{bonnans2000perturbation}, if there exist a pointed closed convex cone $C\subset\bR^m$, a neighborhood $U$ of $x_0$ and a $C^2$ function $\Xi:U\to\bR^m$ such that $\Xi^\prime(x_0):\bR^n\to\bR^m$ is onto, $\Xi(x_0)=0$, and $S\cap U=\{x\in U\mid \Xi(x)\in C\}$. Note that $m\leq n$. If $m=n$, then $\Xi$ is a local $C^2$-diffeomorphism around $x_0$ as $\Xi^\prime(x_0)$ is nonsingular. And if $m<n$, the reduction changes the dimension of space and $\Xi$ is not a local diffeomorphism around $x_0$. However, we show that $\Xi$ can always be extended to a local $C^2$-diffeomorphism in the following lemma.

\begin{lemma}\label{C2hom}
    Let $S\subset\bR^n$ be a nonempty closed convex set which is $C^2$-cone reducible at $\bar{x}\in S$. Then there exist a closed convex cone $C\subset\bR^n$, open neighborhoods $U$ of $\bar{x}$ and $V$ of $0\in\bR^n$ and a $C^2$-diffeomorphism $h:U\to V$ such that $h(\bar{x})=0$ and $h(S\cap U)=C\cap V$.
\end{lemma}
\begin{proof}
    By the definition of $C^2$-cone reduction, there exist a closed convex cone $C_0\subset\bR^m$, an open neighborhood $U_0$ of $\bar{x}$ and a $C^2$ function $\Xi:U_0\to \bR^m$ such that $\Xi'(\bar{x}):\bR^n\to\bR^m$ is onto, $\Xi(\bar{x})=0$ and $S\cap U_0=\{x\in U_0\mid \Xi(x)\in C_0\}$. In particular, $m\leq n$. After a possible relabeling of the coordinates, we may assume without loss of generality that the $m\times m$ submatrix formed by the first $m$ columns of $\Xi^\prime(\bar{x})$ is nonsingular. Define the function $h:U_0\to \bR^n$
    $$
    h(x):=(\Xi_1(x),...,\Xi_m(x),x_{m+1}-\bar{x}_{m+1},...,x_{n}-\bar{x}_n), \quad x\in U_0.
    $$
    So $h'(\bar{x})$ is nonsingular. Since $h$ is $C^2$, there exist open neighborhoods $U\subset U_0$ of $\bar{x}$ and $V$ of $h(\bar{x})=0$ such that $h:U\to V$ is a $C^2$-diffeomorphism. Define the closed convex cone $C\subset\bR^n$ as:
    $$C:=\{y\in\bR^n\mid(y_1,...,y_m)\in C_0\}.$$
    
    We now show that $h(S\cap U)=C\cap V$. For a point $x\in U\subset U_0$, we have $x\in S$ if and only if $\Xi(x)\in C_0$ since $S\cap U_0=\{x\in U_0\mid \Xi(x)\in C_0\}$. So
    \begin{align*}
        y\in h(S\cap U)\iff y\in V\land h^{-1}(y)\in S 
        \iff y\in V\land\Xi(h^{-1}(y))\in C_0\iff y\in C\cap V.
    \end{align*}
    This completes the proof.
\end{proof}

We say that $S\subset\bR^n$ is $C^2$-cone reducible if $S$ is a nonempty closed convex set that is $C^2$-cone reducible at each $x\in S$. According to \cite[Theorem 7.2]{bonnans1998sensitivity} and \cite[Proposition 3.136]{bonnans2000perturbation}, the metric projection onto a $C^2$-cone reducible set is directionally differentiable, which can be combined with the following lemma in our subsequent analysis.
\begin{lemma}[{\cite[Theorem 3.3]{fusek2001isolated}}]\label{isolated}
    Let $f:\bR^n\to\bR^n$ be a Lipschitz continuous and directionally differentiable function, such that $f^{-1}$ has the Aubin property around $(f(x_0),x_0)$. Then $f$ is discrete around $x_0$.
\end{lemma}

Note that with Lemmas \ref{open} and \ref{isolated}, the condition that $f$ is open and discrete around $x_0$ in Lemma \ref{localhom} is easily satisfied, and the most difficult part is $|\ind(f,x_0)|=1$, which is highly related to the structure of $f$ around $x_0$.

\section{Index of \texorpdfstring{$N$}{N}}\label{sec:index}

In this section, we study the function $N$ defined by
$$N(x)=A(x-\Pi_K(x))+B\Pi_K(x),\quad x\in \bR^n,$$
where $K$ is a nonempty closed convex set and $A$, $B$ are $n$ by $n$ matrices. We aim to prove that the index of $N$ at a point $x_0\in K$ is $\pm1$, under the assumption that $N$ is open and discrete around $x_0$. If $x_0\in\ri K$, then $N$ is linear around $x_0$, and hence $|\ind(N,x_0)|=1$ follows directly. {However, when $x_0\in\rb K$, such a direct computation is no longer available. The idea is to use the local constancy of the degree. Namely, if $D$ is a bounded open neighborhood of $x_0$ such that $N^{-1}(N(x_0))\cap\overline D=\{x_0\}$, then $\operatorname{dist}(N(x_0),N(\partial D))>0$. Hence, for all $x$ sufficiently close to $x_0$, the values $N(x)$ and $N(x_0)$ lie in the same connected component of $\mathbb R^n\setminus N(\partial D)$, and therefore
$$
\deg(N|_{\overline D},D,N(x))=\deg(N|_{\overline D},D,N(x_0))=\ind(N,x_0).
$$
However, this still does not make the degree directly computable. The crucial step is to choose $D$ more carefully so that
$$
N^{-1}(N(x))\cap \overline D=\{x\},\quad \forall x\in\{x_0\}\cup(\ri K\cap D).
$$
Then, for every $x\in\ri K\cap D$, $\deg(N|_{\overline D},D,N(x))$ can be computed from the unique preimage $x$, at which $N$ is locally linear, and hence has absolute value one. Taking such an $x\in\ri K\cap D$ sufficiently close to $x_0$ and applying the preceding local constancy argument yields $|\ind(N,x_0)|=1$. Lemma \ref{ANT} provides the key properties for constructing such a domain $D$ in the proof of the main index result, Theorem \ref{indexthm}.}

To prove Lemma \ref{ANT}, we will need the following local structure lemma for the normal cone map of a closed convex set $K$. {The classification of normal vectors recalled in Section \ref{sec:nv} is used here in a precise way: if $u_0\in\rb N_K(x_0)$ and $d=\dim N_K(x_0)>1$, then $T(K,u_0)$ is a proper face of $N_K(x_0)$, so $u_0$ is a $(d-2)$-extreme normal vector. Lemma \ref{extexpthm} then allows us to approximate $u_0$ by $(d-2)$-exposed normal vectors, which yields nearby points with normal cones of dimension strictly smaller than $d$. This lemma provides the key step used in Lemma \ref{ANT} to run a finite descent argument on the dimensions of normal cones.}

\begin{lemma}\label{dimlem}
    Let $K\subset\bR^n$ be a closed convex set with nonempty interior. Then for any $x_0\in \partial K$ and $u_0\in \rb N_K(x_0)$, there exist sequences $\{x_i\}$ and $\{u_i\}$ converging to $x_0$ and $u_0$, respectively, such that $u_i\in N_K(x_i)$ and $\dim N_K(x_i)<\dim N_K(x_0)$ for each $i$.
\end{lemma}
\begin{proof}
    By taking $K:= K \cap \{x\in \bR^n\mid \|x-x_0\|\le r\}$ for some $r>0$ if necessary, we can assume that $K$ is compact.
    Let $d:=\dim N_K(x_0)$. So $d\geq1$ since $x_0\in\partial K$. If $u_0=0$, then we can prove the lemma by choosing $x_i\in\ior K$ and $u_i=0$ and we have $\dim N_K(x_i)=0<d$.
    
    Assume $u_0\neq 0$. Then we have $d\geq2$. The touching cone $T(K,u_0)$ is the smallest face of $N_K(x_0)$ that contains $u_0$. Thus, $\dim T(K,u_0)<d$, since $u_0\in\rb N_K(x_0)$, which implies that $u_0$ is a $(d-2)$-extreme normal vector of $K$. By Lemma \ref{extexpthm}, there is a sequence of $(d-2)$-exposed normal vectors $\{u_k\}$ converging to $u_0$. For each $k$, choose any $\hat{x}_k\in \ri N_K^{-1}(u_k)$. So $\dim N_K(\hat{x}_k)=\dim N_K(\ri N_K^{-1}(u_k))<d$ by the definition of $(d-2)$-exposed normal vectors. Since $K$ is compact and $\hat{x}_k\in K$, there is a convergent subsequence $\{\hat{x}_{k_i}\}$ such that $\hat{x}_{k_i}\to\hat{x}\in K$. The outer semicontinuity of $N_K$ implies $u_0\in N_K(\hat{x})$. If $\dim N_K^{-1}(u_0)=0$, then $\hat{x}=x_0$, which implies that the sequences $\{\hat{x}_{k_i}\}$ and $\{u_{k_i}\}$ satisfy the conditions required in this lemma.
    
    Next, we consider the case that $\dim N_K^{-1}(u_0)>0$. {In this case, the preceding argument may produce points converging to some point in $N_K^{-1}(u_0)$ rather than to the prescribed point $x_0$, so we localize the construction near $x_0$ by introducing $D^\delta$.} Let $L$ be the smallest linear subspace containing $N_K^{-1}(u_0)-x_0$, and let $L^\perp$ be the orthogonal complement of $L$. So $u_0\in N_K(\ri N^{-1}_K(u_0))\subset L^\perp$. For any $\delta>0$, define the closed convex set
    \begin{equation*}
        D^\delta := \{x\in x_0+L\mid\|x-x_0\|\leq\delta\}+L^\perp.
    \end{equation*}
    {Clearly, $x_0\in\mathrm{int}D^\delta$. Since $x_0\in K=\overline{\mathrm{int}K}$, it follows that $\mathrm{int}K\cap\mathrm{int}D^\delta\neq\emptyset$.} So $N_{K\cap D^\delta}(x)=N_K(x)+N_{D^\delta}(x)$ for each $x\in\bR^n$. Note that $u_0\in\rb N_K(x_0)=\rb N_{K\cap D^\delta}(x_0)$. Using the above approach to $x_0$, $u_0$ and $K\cap D^\delta$, there exist $\{x^\delta_k\}$ and $\{u^\delta_k\}$ such that ${x}^\delta_k\to{x}^\delta$, $u^\delta_k\to u_0$, $u^\delta_k\in N_{K\cap D^\delta}({x}^\delta_k)$ and $\dim N_{K\cap D^\delta}({x}^\delta_k)<d$. So $\|{x}^{\delta}-x_0\|\leq\delta$, since ${x}^\delta\in N^{-1}_{K\cap D^\delta}(u_0)=N^{-1}_{K}(u_0)\cap D^\delta\subset (x_0+L)\cap D^\delta$, {where the last inclusion follows from the definition of $L$.}
    
    By utilizing $\{x^\delta_k\}$ and $\{u^\delta_k\}$, we now construct sequences $\{\hat{x}^\delta_k\}$ and $\{\hat{u}^\delta_k\}$ that converge to $x^\delta$ and $u_0$ respectively, while satisfying the conditions $\hat{u}^\delta_k \in N_K(\hat{x}^\delta_k)$ and $\dim N_K(\hat{x}^\delta_k) < d$. If $\|{x}^{\delta}-x_0\|<\delta$, then ${x}^{\delta}\in\ior D^\delta$, which implies that we can find a subsequence $\{{x}^\delta_{k_i}\}\subset\ior D^\delta$. So $N_{K\cap D^\delta}({x}^\delta_{k_i})=N_K({x}^\delta_{k_i})$. Thus, $u_{k_j}^\delta\in N_K({x}^\delta_{k_i})$ and $\dim N_K({x}^\delta_{k_i})<d$. 
    
    Suppose that $\|{x}^{\delta}-x_0\|=\delta$. For each $k$, there exists $v_k\in N_{D^\delta}({x}^\delta_k)$ such that $\hat{u}^\delta_k:=u_k^\delta-v_k\in N_K({x}^\delta_k)$. We claim that $v_k\to0$. Suppose for the sake of contradiction that it is not true. Then there are two cases. One is that there exists a subsequence $\{v_{k_i}\}$ such that $\|v_{k_i}\|\to\infty$. So $\|\hat{u}^\delta_{k_i}\|\to\infty$. Without loss of generality, we can assume that $v_{k_i}\neq0$ and $\hat{u}^\delta_{k_i}\neq0$ for each $i$ and the sequence $\{\hat{u}^\delta_{k_i}/\|\hat{u}^\delta_{k_i}\|\}$ converges to $e$. Since $\hat{u}^\delta_{k_i}\in N_K({x}^\delta_{k_i})$, we have $\langle \hat{u}^\delta_{k_i},{x}^\delta_{k_i}-x_0 \rangle\geq0$ which implies $\langle e,x^\delta-x_0\rangle\geq0$. The outer semicontinuity of $N_{D^\delta}$ implies $v_{k_i}/\|v_{k_i}\|\to({x}^\delta-x_0)/\delta$. So $\langle \hat{u}^\delta_{k_i}/\|\hat{u}^\delta_{k_i}\|,v_{k_i}/\|v_{k_i}\|\rangle\to\langle e,(x^\delta-x_0)/\delta\rangle\geq0$.
    Therefore, $\|u^\delta_{k_i}\|^2=\|\hat{u}^\delta_{k_i}+v_{k_i}\|^2\to\infty$, which is a contradiction. The other case is that there exists a convergent subsequence $v_{k_i}\to v\neq0$. Note that $\langle u_0,{x}^\delta-x_0\rangle=0$, since $u_0\in L^\perp$ and ${x}^\delta-x_0\in L$. By $v_{k_i}/\|v_{k_i}\|\to({x}^\delta-x_0)/\delta$, we have $\langle v,{x}^\delta-x_0\rangle>0$. So $\langle u_0-v,{x}^\delta-x_0\rangle<0$. However, by the outer semicontinuity of $N_K$, we have $u_0-v=\lim \hat{u}^\delta_{k_i}\in N_K({x}^\delta)$, which implies $\langle u_0-v,{x}^\delta-x_0\rangle\geq0$. This is a contradiction. Therefore, $v_k\to0$ and $\hat{u}^\delta_{k}\to u_0$. Moreover, we have $\hat{u}^\delta_{k}\in N_K(x^\delta_k)$ and $\dim N_K({x}^\delta_k)\leq\dim N_{K\cap D^\delta}({x}^\delta_k)<d$.
    
    Therefore, for any $\delta>0$, there exist $\hat{x}^\delta_k\to {x}^\delta$ and $\hat{u}^\delta_k\to u_0$ such that $\|{x}^\delta-x_0\|\leq\delta$, $\hat{u}^\delta_k\in N_K(\hat{x}^\delta_k)$ and $\dim N_K(\hat{x}^\delta_k)<d$. For each $i$, let $\delta_i=1/i$, then we can choose $x_i$ and $u_i$ such that $u_i\in N_K(x_i)$, $\|x_i-x_0\|<2\delta_i$, $\|u_i-u_0\|<\delta_i$ and $\dim N_K(x_i)<d$. So $x_i\to x_0$ and $u_i\to u_0$, which completes the proof.
\end{proof}

The assumption of a nonempty interior in Lemma \ref{dimlem}, which simplifies the proof, can be relaxed. The core argument can be applied to the relative interior of any nonempty closed convex set within its affine hull. This yields a more general result, which we state as a corollary due to its potential independent interest.

\begin{corollary}
    Let $K\subset\bR^n$ be a closed convex set. Then for any $x_0\in \partial K$ and $u_0\in \rb N_K(x_0)$, there exist sequences $\{x_i\}$ and $\{u_i\}$  converging to $x_0$ and $u_0$, respectively, such that $u_i\in N_K(x_i)$ and $\dim N_K(x_i)<\dim N_K(x_0)$ for each $i$.
\end{corollary}
\begin{proof}
    The proof follows by applying Lemma \ref{dimlem} to $K$ within its affine hull, and the details are omitted as it is not essential for the main results of this paper.
\end{proof}

{We next isolate the key geometric argument in the special case where $B=I$ and $K$ has nonempty interior. The resulting lemma will later be applied in the proof of Theorem \ref{indexthm}, after an appropriate transformation and a restriction to the affine hull of $K$. Its proof relies on Lemma \ref{dimlem}, which allows us to pass to nearby boundary points with strictly lower-dimensional normal cones and then argue by finite descent on their dimensions.}
\begin{lemma}\label{ANT}
    Let $K\subset \bR^n$ be a closed convex set with nonempty interior and $A$ be an $n$ by $n$ matrix. Let $N(x)=A(x-\Pi_K(x))+\Pi_K(x)$ for $x\in \bR^n$. Assume that $N$ is open around $x_0\in \partial K$. Then  
    $$AN_K(x_0)\cap\ior T_K(x_0)=\emptyset.$$
\end{lemma}

\begin{proof}
    Let $U$ be an open neighborhood of $x_0$ such that $N|_U$ is an open map. {The proof proceeds by a finite descent on the dimension of the normal cone. To apply Lemma \ref{dimlem}, we first show that whenever the desired conclusion fails at a boundary point, the corresponding normal vector can be chosen from the relative boundary of the normal cone. More precisely, we claim that}
    \begin{equation}\label{rbclaim}
        AN_K(x)\cap\ior T_K(x)\neq\emptyset \ \ \implies\ \ A\,\rb N_K(x)\cap \ior T_K(x)\neq\emptyset,\quad \forall\,x\in \partial K\cap U.
    \end{equation}
    For the sake of contradiction, suppose that there exists a point $\bar{x}\in\partial K\cap U$ such that
    \begin{equation*}
        AN_K(\bar{x})\cap\ior T_K(\bar{x})\neq\emptyset\quad\text{and}\quad A\,\rb N_K(\bar{x})\cap \ior T_K(\bar{x})=\emptyset.
    \end{equation*}
    Let $L$ be the smallest linear subspace containing $N_K(\bar{x})$ and let $P$ be the orthogonal projection matrix onto $L$. Then $y\in\ior T_K(\bar{x})$ if and only if $Py\in\ri(T_K(\bar{x})\cap L)$, since $\ior T_K(\bar{x})=\ri (T_K(\bar{x})\cap L)+L^\perp$. So we have
    \begin{equation}\label{PANK}
        PAN_K(\bar{x})\cap\ri(T_K(\bar{x})\cap L)\neq\emptyset\quad\text{and}\quad PA\,\rb N_K(\bar{x})\cap \ri(T_K(\bar{x})\cap L)=\emptyset.
    \end{equation}
    {Note that $N_K(\bar{x})$ is full-dimensional relative to $L$ by the definition of $L$, and contains no line as $N_K^\circ(\bar{x})=T_K(\bar{x})$ is full-dimensional.} We next show that $PAL=L$. Indeed, if this were false, then $\dim \{u\in L\mid PAu=v\}>0$ for any $v\in PAN_K(\bar{x})\cap\ri(T_K(\bar{x})\cap L)$, thus $\rb N_K(\bar{x})\cap\{u\in L\mid PAu=v\}\neq\emptyset$, i.e., there exists $u^\prime\in\rb N_K(\bar{x})$ such that $PA u^\prime=v$, which contradicts \eqref{PANK}. Thus $PA|_L:L\to L$ is a linear isomorphism. Therefore, $PA\,\rb N_K(\bar{x})=\rb PAN_K(\bar{x})$ and $PAN_K(\bar{x})$ is a closed convex cone and full-dimensional relative to $L$. Combining this with \eqref{PANK}, we have that $(T_K(\bar{x})\cap L)\subset PAN_K(\bar{x})$. Note that 
    $$\dim L=\dim PAN_K(\bar{x})\leq\dim AN_K(\bar{x})\leq \dim N_K(\bar{x})=\dim L.$$
    Thus $\dim PAL=\dim AL=\dim L$. Moreover, $P(AL)=PAL=L$, so $P|_{AL}$ is injective, and hence $AL\cap L^\perp=\{0\}$. In particular, $AN_K(\bar{x})\cap L^\perp=\{0\}$. {So, $AN_K(\bar{x})$ is also a closed convex cone that contains no line, as these properties of $N_K(\bar{x})$ are preserved under the linear isomorphism $A|_L:L\to AL$. 
    Thus $(AN_K(\bar{x}))^\circ$ is full-dimensional. Moreover, $u\in\text{int}[(AN_K(\bar{x}))^\circ]$ if and only if $\langle u,v\rangle<0$ for all $v\in AN_K(\bar{x})\backslash\{0\}$. We claim that $\text{int}[(AN_K(\bar{x}))^\circ]\cap L\neq\emptyset$. Assume on the contrary that $\text{int}[(AN_K(\bar{x}))^\circ]\cap L=\emptyset$. By the separating hyperplane theorem \cite[Theorem 2.39]{rockafellar1998variational}, there exist $u\in\bR^n\backslash\{0\}$ and $\alpha\in\bR$ such that
    $$\text{int}[(AN_K(\bar{x}))^\circ]\subset\{v\in\bR^n\mid\langle u,v\rangle\leq\alpha\}\quad\text{and}\quad L\subset\{v\in\bR^n\mid\langle u,v\rangle\geq\alpha\}.$$
    The conic structure of $\text{int}[(AN_K(\bar{x}))^\circ]$ and the fact that $L$ is a linear subspace force $\alpha=0$. So $u\in ((AN_K(\bar{x}))^\circ)^\circ=AN_K(\bar{x})$ and $u\in L^\perp$, which contradicts $AN_K(\bar{x})\cap L^\perp=\{0\}$. Therefore, there exists $u^\prime\in \text{int}[(AN_K(\bar{x}))^\circ]\cap L$. Define the closed half-space $H:=\{v\in\bR^n\mid\langle u^\prime,v\rangle\leq0\}$. Then $AN_K(\bar{x})\backslash\{0\}\subset \ior H$ and $L^\perp\subset \partial H$.} Thus, we have
    \begin{equation*}
        T_K(\bar{x})=(T_K(\bar{x})\cap L)+L^\perp\subset PAN_K(\bar{x})+L^\perp=AN_K(\bar{x})+L^\perp\subset H.
    \end{equation*}
    Consider the continuous function $f:\bR^n\to\bR^n$ given by $f(x):=A(x-\Pi_K(x))$. Note that $\bar{x}\in f^{-1}(H)$. We claim that $f^{-1}(H)$ is a neighborhood of $\bar{x}$. Suppose that our claim is false. Then there exists a sequence $\{x_i\}$ converging to $\bar{x}$ such that $x_i\notin f^{-1}(H)$ for each $i$. It is clear that $x_i-\Pi_K(x_i)\neq0$ for each $i$ and $\Pi_K(x_i)\to\bar{x}$. Let $u_i:=(x_i-\Pi_K(x_i))/\|x_i-\Pi_K(x_i)\|\in N_K(\Pi_K(x_i))$, and we have $Au_i=f(x_i)/\|x_i-\Pi_K(x_i)\|\notin H$ for each $i$, since $H$ is a cone. Passing to a subsequence if necessary, assume that $u_i\to\bar{u}$. So $A\bar{u}\notin\ior H$. But the outer semicontinuity of $N_K$ implies $A\bar{u}\in AN_K(\bar{x})\backslash\{0\}\subset\ior H$, where $A\bar{u}\neq0$ since $\bar{u}\in L\backslash\{0\}$ and $\dim AL=\dim L$. Therefore, $f^{-1}(H)$ must be a neighborhood of $\bar{x}$ and there exists an open neighborhood $V\subset f^{-1}(H)\cap U$ of $\bar{x}$. We have $\bar{x}=N(\bar{x})\in N(V)$ and
    \begin{equation*}
        N(V)= \bigcup_{x\in V}(\Pi_K(x)+f(x))\subset K+H\subset \bar{x}+T_K(\bar{x})+H=\bar{x} + H.
    \end{equation*}
    It is clear that $\bar{x}$ is not an interior point of $N(V)$. However, since $N|_U$ is an open map, $N(V)$ is an open neighborhood of $\bar{x}$. This is a contradiction. So claim \eqref{rbclaim} is proved.

    We now combine \eqref{rbclaim} with Lemma \ref{dimlem}. Suppose, to the contrary, that
    \begin{equation*}
        AN_K(x_0)\cap\ior T_K(x_0)\neq\emptyset.
    \end{equation*}
    Then $A\,\rb N_K(x_0)\cap\ior T_K(x_0)\neq\emptyset$.
    Let $d:=\dim N_K(x_0)\geq 1$. If $d=1$, we have $\rb N_K(x_0)=\{0\}$ and $A\,\rb N_K(x_0)\cap\ior T_K(x_0)=\emptyset$, which is a contradiction. Assume $d>1$. Then there exists a normal vector $u_0\in\rb N_K(x_0)\backslash\{0\}$ such that $Au_0\in\ior T_K(x_0)$.
    By the inner semicontinuity of the convex-valued map $T_K$ and \cite[Theorem 5.9]{rockafellar1998variational}, there is an open neighborhood $W\subset\bR^n\times\bR^n$ of $(x_0,u_0)$, such that $Au\in \ior T_K(x)$ for any $(x,u)\in W$. So for any $(x,u)\in \gph N_K$ sufficiently close to $(x_0,u_0)$, we have $AN_K(x)\cap\ior T_K(x)\neq\emptyset$.
    Therefore, by Lemma \ref{dimlem}, we can find $x^\prime_0\in \partial K\cap U$ such that 
    \begin{equation*}
        AN_K(x^\prime_0)\cap\ior T_K(x^\prime_0)\neq\emptyset\quad\text{and}\quad\dim N_K(x^\prime_0)<\dim N_K(x_0).
    \end{equation*}

    We may now repeat the same argument with $x_0^\prime$ in place of $x_0$. At each repetition, the newly obtained point has a normal cone whose dimension is strictly smaller than that at the preceding point, so the process terminates after finitely many steps. Hence we obtain $x^{\prime\prime}_0\in\partial K\cap U$ such that $AN_K(x^{\prime\prime}_0)\cap\ior T_K(x^{\prime\prime}_0)\neq\emptyset$ and $\dim N_K(x^{\prime\prime}_0)=1$, which contradicts the fact that the case $d=1$ has already been ruled out. Therefore, we have $AN_K(x_0)\cap \ior T_K(x_0)=\emptyset$.
\end{proof}

Under the assumption of Lemma \ref{ANT}, there exists an open neighborhood $U$ of $x_0$ such that $N$ is open around each $x\in\partial K\cap U$. So for each $x\in \partial K\cap U$, we have $AN_K(x)\cap \ior T_K(x)=\emptyset$, which implies $(x+N_K(x))\cap N^{-1}(\ior K)=\emptyset$. Thus, for any $x\in \ri K\cap U$, $N^{-1}(N(x))=\{x\}$. If $N$ is discrete at $x_0$, we can prove $|\ind(N,x_0)|=1$ by the local constancy of the degree. {Motivated by this observation, we establish the following index theorem for $N$ in a more general form, allowing both an additional matrix $B$ and a closed convex set $K$ whose interior may be empty.}

\begin{theorem}\label{indexthm}
    Let $K\subset \bR^n$ be a nonempty closed convex set and $A,B$ be $n$ by $n$ matrices. Let 
    $$N(x)=A(x-\Pi_K(x))+B\Pi_K(x), \quad x\in\bR^n.$$
    Assume that $N$ is open around $x_0\in K$ and discrete at $x_0$. Then $|\ind(N,x_0)|=1$.
\end{theorem}

\begin{proof}
    Let $N_0(x):=N(x+x_0)-N(x_0)=A(x-\Pi_{K-x_0}(x))+B\Pi_{K-x_0}(x)$ for $x\in \bR^n$. Then $N_0$ is open around $0\in K-x_0$ and discrete at $0$. Moreover, $\ind(N,x_0)=\ind(N_0,0)$. Thus, without loss of generality, we assume that $x_0=0$. Then $0\in K$ and $N(0)=0$.
    
    Let $L$ be the smallest linear subspace containing $K$. Since $0\in K$, this is also the affine hull of $K$. Also, let $P$ be the orthogonal projection matrix onto $L$. By assumption, there exists a bounded open neighborhood $U$ of $0$ such that $N|_U$ is an open map and $N^{-1}(0)\cap U=\{0\}$. Since $0\in K$, we know that $\ri K\cap U\neq\emptyset$. For any $y\in\ri K$, we have $N_K(y)=L^\perp$. Hence $\Pi_K(x)=Px$ and $N(x)=(A-AP+BP)x$ for any $x\in \Pi_K^{-1}(\ri K)$. Moreover,
    \[
        \Pi_K^{-1}(\ri K)
        = \bigcup_{y\in \ri K}\Pi_K^{-1}(y)
        = {\bigcup_{y\in \ri K}(I+N_K)(y)}
        = \bigcup_{y\in \ri K}y + L^\perp
        = \ri K + L^\perp.
    \]
    {The relative openness of $\ri K$ in $L$, together with the orthogonal decomposition of $\bR^n$ into $L$ and $L^\perp$, implies that $\ri K+L^\perp$ is open in $\bR^n$.} Since the function $N$ is linear on the open set $\Pi_K^{-1}(\ri K)$ and is open around any $x\in \ri K\cap U\subset \Pi_K^{-1}(\ri K)$, we obtain that the matrix $A-AP+BP$ is nonsingular. If $0\in\ri K$, then 
    $$
    |\ind(N,0)|=|\deg(N|_{\overline{U}},U,0)|=|\sgn\det N^\prime(0)|=1.
    $$ 
    It remains to consider the case $0\in\rb K$. Let $M=(A-AP+BP)^{-1}$. Then
    \begin{equation}\label{eq:simple}
    MAL^\perp= L^\perp \ \ {\rm and} \ \ MBx=M[(Ax-APx) +BPx] =x, \ \  \forall\,  x\in L.
    \end{equation}
    Since $K\subset L$, we have $\Pi_K(x)=\Pi_K(Px)$ for all $x\in\bR^n$. Thus, for any $x\in \bR^n$, 
    $$
    PMN(x)=PMA(x-\Pi_K(x))+\Pi_K(x)=PMA(Px-\Pi_K(Px))+\Pi_K(Px).
    $$
    Let $\widehat{A}:L\to L$ be the linear map given by $\widehat{A}z=PMAz$, $z\in L$. We now reduce the problem to the subspace $L$. Define the function $\widehat{N}:L\to L$ by
    $$
    \widehat{N}(x)=PMN(x)=\widehat{A}(x-\Pi_K(x))+\Pi_K(x), \quad x\in L.
    $$
    For each open set $V\subset U\cap L$ in the linear subspace $L$, we have $\widehat{N}(V)=PMN((V+L^\perp)\cap U)$, which is open in $L$, since $N|_U:U\to\bR^n$ and $P:\bR^n\to L$ are open maps. Therefore, $\widehat{N}|_{U\cap L}$ is an open map. Since the interior of $K$ relative to $L$, namely $\ri K$, is nonempty, Lemma \ref{ANT} applied to the function $\widehat{N}:L\to L$ gives that for each $x\in \rb K\cap U$,
    \begin{equation}
    \label{eq:ANT}
        \widehat{A}N_K^L(x)\cap \mathrm{int}^L\, T^L_K(x)=\emptyset,
    \end{equation}
    where $N^L_K(x)$ and $T^L_K(x)$ are the normal cone and the tangent cone to $K$ at $x$ relative to $L$, respectively, and $\mathrm{int}^L$ means the interior relative to $L$. Moreover, $N^L_K(x)=PN_K(x)$ and $\mathrm{int}^L\, T^L_K(x)=\ri T_K(x)$. We now translate the relative conclusion \eqref{eq:ANT} back to the original setting. Fix $x\in\rb K\cap U$. Then we have
    $$
    PMAPN_K(x)\cap \ri T_K(x)=\emptyset.
    $$
    Taking the Minkowski sum with $L^\perp$ preserves the empty intersection, since both original sets lie in $L$, that is,
    $$
    (PMAPN_K(x) + L^\perp)\cap (\ri T_K(x)+L^\perp)=\emptyset.
    $$
    Thus, by using \eqref{eq:simple}, we have 
    $$
    PMAPN_K(x)+L^\perp = MAPN_K(x)+L^\perp= MA(PN_K(x)+L^\perp)=MAN_K(x),
    $$
    which further implies that $MAN_K(x)\cap(\ri T_K(x)+L^\perp)=\emptyset$. Applying $M^{-1}$ to both sets, and using the identities $M^{-1}\ri T_K(x) = B\,\ri T_K(x)$ and $M^{-1}L^\perp = AL^\perp$ given by \eqref{eq:simple}, we obtain
    $$
    AN_K(x)\cap(B\,\ri T_K(x)+AL^\perp)=\emptyset.
    $$
    Consequently, for any $x \in \rb K \cap U$,
    \begin{equation}
        \label{eq:simple2}N(x+N_K(x))\cap N(\ri K)\subset (Bx+AN_K(x))\cap(B\,\ri K+AL^\perp)=\emptyset,
    \end{equation} 
    where we use the fact that $\ri K-x \subset \ri T_K(x)$. 
    
    Let $U^\prime:=U\cap \Pi_K^{-1}(U)$. Then $U^\prime$ satisfies all the requirements of $U$, and $\Pi_K(U^\prime)\subset U^\prime$ since $\Pi_K(\Pi_K(U^\prime))=\Pi_K(U^\prime)\subset U$. By replacing $U$ with $U^\prime$, we can assume that $U\subset\Pi_K^{-1}(U)$. We claim that for each $z\in \ri K\cap U$,
    $$
    N^{-1}(N(z))\cap U=\{z\}.
    $$
    To see this, fix $z\in \ri K\cap U$ and suppose that $w\in N^{-1}(N(z))\cap U$. Then there are two cases: $\Pi_K(w)\in\ri K$ and $\Pi_K(w)\in\rb K$. If $\Pi_K(w)\in\ri K$, then we have that $M^{-1}w=N(w)=N(z)=M^{-1}z$, which implies $w = z$.  If $\Pi_K(w)\in\rb K$, then $\Pi_K(w)\in \rb K\cap U$ as $w\in U\subset\Pi_K^{-1}(U)$. Therefore, by using \eqref{eq:simple2}, we know that $$N(\Pi_K(w)+N_K(\Pi_K(w)))\cap N(\ri K)=\emptyset,$$ which, together with the fact that $ w \in \Pi_K(w) + N_K(\Pi_K(w))$, 
    implies $N(w)\notin N(\ri K)$. However, this contradicts our assumption that $N(w) = N(z) \in N(\ri K)$. Thus, $w$ must be equal to $z$. Therefore, $N^{-1}(N(z))\cap U=\{z\}$ for any $z\in \ri K\cap U$.
    
    Let $D$ be a bounded open neighborhood of $0$ with $\overline{D} \subset U$. For any $z \in \ri K \cap D$, since $N^{-1}(N(z))\cap U=\{z\}$, we have $N(z)\notin N(\partial D)$ and 
    $$
    \deg(N|_{\overline{D}},D,N(z))=\sgn\det N'(z)=\sgn\det(A-AP+BP).
    $$
    It follows from $N^{-1}(0) \cap U = \{0\}$ that $0 \notin N(\partial D)$. The image $N(\partial D)$ is compact because it is the continuous image of the compact set $\partial D$. Thus, $\dist(0, N(\partial D)) > 0$. Since there exists a sequence $\{z_i\}\subset\ri K\cap D$ converging to $0\in\rb K\cap D$, the local constancy of the degree gives
    $$
    |\ind(N,0)|=|\deg(N|_{\overline{D}},D,N(0))|=|\sgn\det(A-AP+BP)|=1,
    $$
    which completes the proof.
\end{proof}

\section{Equivalence of the Aubin property and strong regularity}\label{sec:eqv}

In this section, we establish the main results of the paper, proving the equivalence between the Aubin property and strong regularity for generalized equations over $C^2$-cone reducible sets. Our proof strategy proceeds in three main steps. First, in Theorem \ref{NABhom}, we tackle a canonical form of the function $N(x)=A(x-\Pi_S(x))+B\Pi_S(x)$, where the core of our degree-theoretic argument, Theorem \ref{indexthm}, is applied. Second, we extend this result in Theorem \ref{NFGhom} to a more general form in which the constant matrices $A$ and $B$ are replaced by $C^1$ functions. Finally, in Theorem \ref{main}, we show how the original generalized equation can be transformed to fit the structure of Theorem \ref{NFGhom}, thus completing the proof of our main claim.

Although Theorem \ref{indexthm} applies to an arbitrary closed convex set, the reference point must lie in that set. In Theorem \ref{NABhom}, however, the local analysis is centered at the point $x_0$, which is not necessarily contained in $S$. Therefore, in Theorem \ref{NABhom}, we first need to use the $C^2$-cone reducibility of $S$ to transform this reference point. Lemma \ref{C2hom} provides a local $C^2$-diffeomorphism $h$ between $S$ near $\Pi_S(x_0)$ and a closed convex cone $C$ near the origin. However, this diffeomorphism acts only on points of $S$ and does not directly transform $\Pi_S$ appearing in $N$ for points near $x_0$. The following lemma bridges this gap. It lifts the diffeomorphism $h$ near $\Pi_S(x_0)$ to a homeomorphism $g$ near $x_0$ such that
$$
\Pi_C\circ g=h\circ \Pi_S.
$$
Consequently, under the change of variables $y=g(x)$, $\Pi_S$ is transferred to $h^{-1}\circ \Pi_C$, and the reference point is transformed to $y_0:=g(x_0)$. Moreover, the above relation gives $\Pi_C(y_0)=0$, i.e., $y_0\in C^\circ$. Although the change of variables $y=g(x)$ provided in the following lemma introduces nonlinearities associated with $h$, it is the key first step in the proof of Theorem \ref{NABhom}.

{
\begin{lemma}
\label{lem:lift}
    Let $S,C\subset \bR^n$ be closed convex sets. Suppose that there exist open sets $U,V\subset \bR^n$ and a $C^1$-diffeomorphism $h:U\to V$ such that
    $$h(S\cap U)=C\cap V\neq\emptyset.$$
    Then there exists a homeomorphism $g:\Pi^{-1}_S(U)\to\Pi^{-1}_C(V)$ satisfying
    $$\Pi_C(g(x))=h(\Pi_S(x)),\quad \forall x\in\Pi^{-1}_S(U).$$
    More precisely, $g$ is defined by
    \begin{equation}
    \label{g}
        g(x):=[\nabla h(\Pi_S(x))]^{-1}(x-\Pi_S(x))+h(\Pi_S(x)), \quad x\in \Pi_S^{-1}(U).
    \end{equation}
    The inverse of $g$ is given by
    $$g^{-1}(y)=\nabla h(h^{-1}(\Pi_C(y)))(y-\Pi_C(y))+h^{-1}(\Pi_C(y)), \quad y\in \Pi^{-1}_C(V).$$
\end{lemma}

\begin{proof}
    Define $g:\Pi_S^{-1}(U)\to W$ by \eqref{g}, where $W:=g(\Pi_S^{-1}(U))\subset\bR^n$. By \cite[Exercise 6.7]{rockafellar1998variational}, for every $z\in S\cap U$, $N_S(z)=\nabla h(z)N_C(h(z))$. For any $x\in\Pi^{-1}_S(U)$, since $x-\Pi_S(x)\in N_S(\Pi_S(x))$ and $\Pi_S(x)\in S\cap U$, we have
    $$
    [\nabla h(\Pi_S(x))]^{-1}(x-\Pi_S(x))\in[\nabla h(\Pi_S(x))]^{-1} N_S(\Pi_S(x))=N_C(h(\Pi_S(x))).
    $$
    Thus, by the formula of $g$ in \eqref{g}, 
    $$
    \Pi_C(g(x))=h(\Pi_S(x)),\quad \forall x\in\Pi^{-1}_S(U).
    $$
    In particular, $W\subset \Pi_C^{-1}(V)$. Conversely, let $y_0\in \Pi_C^{-1}(V)$. Set $z_0:=h^{-1}(\Pi_C(y_0))$.
    So $z_0\in S\cap U$. Moreover, $\nabla h(z_0)(y_0-\Pi_C(y_0))\in N_S(z_0)$. Define $x_0:=\nabla h(z_0)(y_0-\Pi_C(y_0))+z_0$. Then $\Pi_S(x_0)=z_0\in U$. Hence $x_0\in \Pi_S^{-1}(U)$, and
    $$g(x_0)=[\nabla h(z_0)]^{-1}(x_0-z_0)+h(z_0)=y_0.$$
    Therefore, $y_0\in W$, which implies $W=\Pi_C^{-1}(V)$.

    Finally, let $y\in \Pi_C^{-1}(V)$ and write $y=g(x)$. Then $\Pi_S(x)=h^{-1}(\Pi_C(y))$. Substituting this into \eqref{g} and solving for $x$ yields
    $$
    g^{-1}(y)=\nabla h(h^{-1}(\Pi_C(y)))(y-\Pi_C(y))+h^{-1}(\Pi_C(y)), \quad y\in \Pi_C^{-1}(V).
    $$
    Hence $g$ is a homeomorphism between $\Pi^{-1}_S(U)$ and $\Pi^{-1}_C(V)$.
\end{proof}
}

\begin{figure}
    \centering
    \caption{A schematic illustration of the lift in Lemma \ref{lem:lift}}
    \label{fig:lift}
    \usetikzlibrary{arrows.meta}
    \usetikzlibrary{calc}
    \usetikzlibrary{patterns}
    \vspace{3mm}
    \begin{tikzpicture}[>=Latex,scale=1.0]
    
    \coordinate (TL) at (0,5);   
    \coordinate (TR) at (8,5);   
    \coordinate (BL) at (0,0);   
    \coordinate (BR) at (8,0);   
    
    \def\Xmin{-1.7}
    \def\Xmax{ 1.7}
    \def\Ymin{-2}
    \def\Ymax{ 2}
    
    \def\rU{0.9}
    \def\Tnormal{3} 
    
    \def\paray#1{(#1*#1)}
    
    \begin{scope}[shift={(TL)}]
      \clip (\Xmin,\Ymin) rectangle (\Xmax,\Ymax);
    
      \draw[very thick] plot[domain=\Xmin:\Xmax,samples=90] (\x,{\paray{\x}});
    
      \path[pattern=north west lines,pattern color=gray!70]
        plot[domain=\Xmin:\Xmax,samples=90] (\x,{\paray{\x}})
        -- (\Xmax,\Ymax) -- (\Xmin,\Ymax) -- cycle;
      \node at (-0.9,1.7) {$S$};
    
      \path[pattern=north east lines,pattern color=blue!70] (0,0) circle (\rU);
      \node[blue!70!black] at (0,1.1) {$U$};
      \fill (0,0) circle (2pt) node[below left] {$\Pi_S(x_0)$};
    
    \end{scope}
    
    \begin{scope}[shift={(TR)}]
      \clip (\Xmin,\Ymin) rectangle (\Xmax,\Ymax);
    
      \pgfmathsetmacro{\xr}{sqrt(\Ymax)}
      \pgfmathsetmacro{\xl}{-sqrt(\Ymax)}
      \draw[very thick] (\xl,0)--(\xr,0);
      \node at (0.9,1.7) {$C$};
    
      \begin{scope}
        \clip (\xl,0) -- (\xr,0) -- (\xr,\Ymax) -- (\xl,\Ymax) -- cycle;
        \path[pattern=north west lines,pattern color=gray!70]
          (\Xmin,\Ymin) rectangle (\Xmax,\Ymax);
      \end{scope}
    
      \path[pattern=north east lines,pattern color=blue!70]
        plot[domain=0:360,samples=80,variable=\t] 
          ({\rU*cos(\t)},{\rU*sin(\t)-(\rU*cos(\t))^2}) -- cycle;
      \node[blue!70!black] at (0,1.1) {$V$};
      \fill (0,0) circle (2pt) node[below right] {$\Pi_C(y_0)$};
    
    \end{scope}
    
    \begin{scope}[shift={(BL)}]
      \clip (\Xmin,\Ymin) rectangle (\Xmax,\Ymax);
    
      \draw[very thick] plot[domain=\Xmin:\Xmax,samples=90] (\x,{\paray{\x}});
    
    
      \begin{scope}
        \clip (0,0) circle (\rU);
        \path[fill=orange!35,opacity=0.7]
          plot[domain=\Xmin:\Xmax,samples=100] (\x,{\paray{\x}})
          -- (\Xmax,\Ymax) -- (\Xmin,\Ymax) -- cycle;
      \end{scope}
    
      \def\xw{0.73}
      \path[fill=orange!35,opacity=0.7]
        plot[domain=-\xw:\xw,samples=90] (\x,{\paray{\x}})
        --
        plot[domain=\xw:-\xw,samples=90]
          ({\x + 2*\x*\Tnormal},{\paray{\x} - \Tnormal})
        -- cycle;
    
      \foreach \xx in {-0.6,-0.4,-0.2,0,0.2,0.4,0.6}{
        \draw[orange!70!black,opacity=0.9]
          ({\xx},{\paray{\xx}}) -- ({\xx+2*\xx*\Tnormal},{\paray{\xx}-\Tnormal});
      }
    
      \coordinate (x0L) at (0,-1.0);
      \coordinate (pKL)  at (0,0);
      \fill (x0L) circle (2pt) node[below left] {$x_0$};
      \fill (pKL) circle (2pt) node[below left] {$\Pi_S(x_0)$};
    
    \end{scope}
    
    \begin{scope}[shift={(BR)}]
      \clip (\Xmin,\Ymin) rectangle (\Xmax,\Ymax);
    
      \pgfmathsetmacro{\xr}{sqrt(\Ymax)}
      \pgfmathsetmacro{\xl}{-sqrt(\Ymax)}
      \draw[very thick] (\xl,0)--(\xr,0);
    
    
      \begin{scope}
        \clip
          plot[domain=0:360,samples=80,variable=\t]
            ({\rU*cos(\t)},{\rU*sin(\t)-(\rU*cos(\t))^2}) -- cycle;
        \path[fill=orange!35,opacity=0.7] (\Xmin,0) rectangle (\Xmax,\Ymax);
      \end{scope}
    
      \def\wCurt{0.73}
      \path[fill=orange!35,opacity=0.7] (-\wCurt,\Ymin) rectangle (\wCurt,0);
    
      \foreach \xx in {-0.6,-0.4,-0.2,0,0.2,0.4,0.6}{
        \draw[orange!70!black,opacity=0.9]
          (\xx,0) -- (\xx,\Ymin);
      }
    
      \coordinate (y0R) at (0,-1.0);
      \coordinate (pCR)  at (0,0);
      \fill (y0R) circle (2pt) node[below right] {$y_0$};
      \fill (pCR) circle (2pt) node[below right] {$\Pi_C(y_0)$};
    
    \end{scope}
    
    
    \draw[<->,blue!70!black,very thick]
      ($(TL)+(0,0.6)$) .. controls ($(TL)+(2.0,1.4)$) and ($(TR)+(-2.0,1.4)$) ..
      ($(TR)+(0,0.6)$)
      node[midway,above,yshift=0.5ex] {$h$}
      node[midway,below,yshift=-0.5ex] {$h(S\cap U)=C\cap V$};
    
    \draw[<->,orange!70!black,very thick]
      ($(BL)+(0,-1.6)$) .. controls ($(BL)+(2.0,-2.2)$) and ($(BR)+(-2.0,-2.2)$) ..
      ($(BR)+(0,-1.6)$)
      node[midway,below,yshift=-0.4ex] {$g$}
      node[midway,above,yshift=0.5ex] {$g(x_0)=y_0$};
    
    \node[orange!70!black] at (0,-2.3) {$\Pi_S^{-1}(U)$};
    \node[orange!70!black] at (8,-2.3) {$\Pi_C^{-1}(V)$};
    
    
    \draw[->,very thick]
      ($(BL)+(0,1)$) -- ($(TL)+(0,-0.2)$) node[midway,left] {$\Pi_S$};
    \draw[->,very thick]
      ($(BR)+(0,1)$) -- ($(TR)+(0,-0.2)$) node[midway,right] {$\Pi_C$};
    
    \node[overlay] at (current bounding box.center) {$\Pi_C(g(x))=h(\Pi_S(x))$};
    
    \end{tikzpicture}
\end{figure}

{
The geometric meaning of Lemma \ref{lem:lift} is illustrated in Figure \ref{fig:lift}. We now turn to the proof of Theorem \ref{NABhom}. For a $C^2$-cone reducible set $S$, Lemme \ref{C2hom} provides a $C^2$-diffeomorphism $h$ reducing $S$ to a cone $C$ around $\Pi_S(x_0)$. Lemma \ref{lem:lift} then lifts $h$ to a homeomorphism $g$ around $x_0$ so that $\Pi_C\circ g=h\circ\Pi_S$. In fact, we can use $y=g(x)$ to change variables around $x_0$ and study $N_0=N\circ g^{-1}$ around $y_0=g(x_0)$ instead, without changing the problem.

Using $\Pi_C\circ g=h\circ\Pi_S$ and Moreau's decomposition, $N_0$ can be written in terms of $\Pi_{C^\circ}$, with $y_0\in C^\circ$. Thus the problem is brought into the framework of Theorem \ref{indexthm}, except for the nonlinear terms introduced by $h$. The rest of the proof removes the nonlinearities through three strictly stationary perturbations at $y_0$. By Lemmas \ref{invmapthm} and \ref{keepindex}, the index is preserved throughout and the final form is covered by Theorem \ref{indexthm}, we obtain $|\ind(N_0,y_0)|=1$, which yields $N$ is a local homeomorphism at $x_0$.
}

\begin{theorem}\label{NABhom}
    Let $S\subset\bR^n$ be a $C^2$-cone reducible set and $A,B$ be $n$ by $n$ matrices. Let 
    $$N(x)=A(x-\Pi_S(x))+B\Pi_S(x), \quad x\in \bR^n.$$
    Assume that $N^{-1}$ has the Aubin property around $(N(x_0),x_0)$. Then $N$ is a local homeomorphism at $x_0$.
\end{theorem}
\begin{proof}
    By \cite[Theorem 7.2]{bonnans1998sensitivity}, $\Pi_S$ is directionally differentiable since $S$ is $C^2$-cone reducible. Consequently, $N$ is also directionally differentiable. Moreover, $N$ is Lipschitz continuous. Thus, by Lemma \ref{isolated}, $N$ is discrete around $x_0$.
    
    By Lemma \ref{C2hom}, there exist a closed convex cone $C\subset\bR^n$, open neighborhoods $U\subset\bR^n$ of $\Pi_S(x_0)$ and $V\subset\bR^n$ of $0$ and a $C^2$-diffeomorphism $h:U\to V$ such that $h(\Pi_S(x_0))=0$ and $h(S\cap U)=C\cap V$. Then, there is a homeomorphism $g:\Pi_S^{-1}(U)\to \Pi_C^{-1}(V)$ associated with $h$ defined by \eqref{g} in Lemma \ref{lem:lift} such that
    $$
    \Pi_S(g^{-1}(y))=h^{-1}(\Pi_C(y)),\quad\forall y\in \Pi_C^{-1}(V),
    $$
    and the inverse $g^{-1}$ is given by
    $$
    g^{-1}(y)=\nabla h(h^{-1}(\Pi_C(y)))(y-\Pi_C(y))+h^{-1}(\Pi_C(y)), \quad y\in \Pi^{-1}_C(V).
    $$
    
    In what follows, the variable $y$ is always taken from $\Pi_C^{-1}(V)$. Define $N_0$ by
    $$
    N_0(y):=N(g^{-1}(y))=A\nabla h(h^{-1}(\Pi_C(y)))(y-\Pi_C(y))+Bh^{-1}(\Pi_C(y)).
    $$
    Let $y_0:=g(x_0)$. Since $h$ is a $C^2$-diffeomorphism, $h$, $h^{-1}$, $\nabla h$ and $[\nabla h]^{-1}$ are continuously differentiable. Hence, from the formulas for $g$ and $g^{-1}$ in Lemma \ref{lem:lift}, $g$ is Lipschitz continuous around $x_0$ and $g^{-1}$ is Lipschitz continuous around $y_0$. Therefore, from the definitions, we can directly check that $N^{-1}$ has the Aubin property around $(N(x_0),x_0)$ if and only if $N_0^{-1}$ has the Aubin property around $(N_0(y_0),y_0)$; $N$ is a local homeomorphism at $x_0$ if and only if $N_0$ is a local homeomorphism at $y_0$; and $N_0$ is discrete around $y_0$ since $N$ is discrete around $x_0$. Moreover, $\Pi_C(y_0)=h(\Pi_S(x_0))=0$, and hence $y_0\in C^\circ$. Using Moreau's decomposition $\Pi_C +\Pi_{C^\circ}=I$, we can rewrite $N_0$ as
    $$
    N_0(y)=Bh^{-1}(y-\Pi_{C^\circ}(y))+A\nabla h(h^{-1}(y-\Pi_{C^\circ}(y)))\Pi_{C^\circ}(y).
    $$
    
    {We next transform $N_0$, by adding strictly stationary perturbations at $y_0$, into the form required in Theorem \ref{indexthm}.} Since $h^{-1}$ is continuously differentiable around $0$, the following function
    $$
    \psi_1(y):=B[(h^{-1})'(0)(y-\Pi_{C^\circ}(y))-h^{-1}(y-\Pi_{C^\circ}(y))]
    $$
    is strictly stationary at $y_0$ by Lemma \ref{ss}. Let $\widehat{B}:=B(h^{-1})'(0)$ and
    $$
    N_1(y):=N_0(y)+\psi_1(y)=\widehat{B}(y-\Pi_{C^\circ}(y))+A\nabla h(h^{-1}(y-\Pi_{C^\circ}(y)))\Pi_{C^\circ}(y).
    $$
    It is clear that $\nabla h\circ h^{-1}$ is continuously differentiable around $0$. So we can define a function $\psi_2$ by
    $$
    \psi_2(y):=A[(\nabla h\circ h^{-1})^\prime(0)(y-\Pi_{C^\circ}(y ))+\nabla h(h^{-1}(0))-\nabla h(h^{-1}(y-\Pi_{C^\circ}(y)))]\Pi_{C^\circ}(y).
    $$
    It follows from Lemma \ref{ss} that $\psi_2$
    is strictly stationary at $y_0$.  Let $M:\bR^n\to\bR^{n\times n}$ be a linear map such that for any $a,b\in\bR^n$, $M(a)b:=A(\nabla h\circ h^{-1})^\prime(0)(b)a$, and let $\widehat{A}:=A\nabla h(h^{-1}(0))$. Define
    $$
    N_2(y):=N_1(y)+\psi_2(y)=(\widehat{B}+M(\Pi_{C^\circ}(y )))(y-\Pi_{C^\circ}(y))+\widehat{A}\Pi_{C^\circ}( y).
    $$
    Consider the function
    $$
    \psi_3(y):=(M(y_0)-M(\Pi_{C^\circ}(y)))(y-\Pi_{C^\circ}(y)).
    $$
    According to Lemma \ref{ss}, we have that $\psi_3$ is strictly stationary at $y_0$. Define
    $$
    N_3(y):=N_2(y)+\psi_3(y)=(\widehat{B}+M(y_0))(y-\Pi_{C^\circ}(y))+\widehat{A}\Pi_{C^\circ}(y).
    $$
    This is precisely the form covered by Theorem \ref{indexthm}.
    
    As $\psi_1$, $\psi_2$ and $\psi_3$ are strictly stationary at $y_0$, Lemma \ref{ss} implies that the difference $N_3-N_0=\psi_1+\psi_2+\psi_3$ is strictly stationary at $y_0$. Thus, by Lemmas \ref{invmapthm} and \ref{keepindex}, $N_3^{-1}$ has the Aubin property around $(N_3(y_0),y_0)$, $N_3$ is discrete at $y_0$, and $\ind(N_0,y_0)=\ind(N_3,y_0)$. Lemma \ref{open} implies that both $N_0$ and $N_3$ are open around $y_0$. Applying Theorem \ref{indexthm} to $N_3$ around $y_0$, we obtain $|\ind(N_3,y_0)|=1$. Hence $|\ind(N_0,y_0)|=1$. By Lemma \ref{localhom}, $N_0$ is a local homeomorphism at $y_0$, and therefore $N$ is a local homeomorphism at $x_0$.
\end{proof}

By Lemma \ref{ss}, we can easily extend Theorem \ref{NABhom} to a more general case.

\begin{theorem}\label{NFGhom}
    Let $S\subset\bR^n$ be a $C^2$-cone reducible set and $f,g:\bR^n\to\bR^n$ be $C^1$ functions. Let 
    $$N(x)=f(x-\Pi_S(x))+g(\Pi_S(x)), \quad x\in \bR^n.$$
    Assume that  $N^{-1}$ has the Aubin property around $(N(x_0),x_0)$. Then $N$ is a local homeomorphism at $x_0$.
\end{theorem}
\begin{proof}
    Define $N_0:\bR^n\to\bR^n$ by 
    $$N_0(x):=f'(x_0-\Pi_S(x_0))(x-\Pi_S(x))+g'(\Pi_S(x_0))\Pi_S(x),\quad x\in\bR^n.$$
    By Lemma \ref{ss}, $N_0-N$ is strictly stationary at $x_0$. The proof is completed by applying Lemma \ref{invmapthm} and Theorem \ref{NABhom}.
\end{proof}

We are now fully equipped to prove the main theorem of this paper. The following proof demonstrates that the solution map $\Phi^{-1}$ for the generalized equation can be related, through a homeomorphism, to a function of the form analyzed in Theorem \ref{NFGhom}. This allows us to translate the Aubin property of $\Phi^{-1}$ into the local homeomorphism of the related function, which in turn implies the strong regularity of $\Phi^{-1}$.

\begin{theorem}\label{main}
    Let $S\subset\bR^n$ be a $C^2$-cone reducible set, and let $\varphi:\bR^n\to\bR^n$ be a $C^1$ function. Let $\Phi:\bR^n\rightrightarrows\bR^n$ be a set-valued map in one of the following forms:
    \begin{equation*}
        \Phi(x)=\varphi(x)+N_S(x)   
        \quad\mathrm{or}\quad
        \Phi(x)=\varphi(x)+N^{-1}_S(x), \quad x \in \bR^n. 
    \end{equation*}
    Then the Aubin property and strong regularity for $\Phi^{-1}$ around $(y_0,x_0)\in\gph\Phi^{-1}$ are equivalent.
\end{theorem}
\begin{proof}
    The implication from strong regularity to the Aubin property is immediate. Conversely, assume that $\Phi^{-1}$ has the Aubin property around $(y_0,x_0)\in\gph\Phi^{-1}$. It remains to show that $\Phi^{-1}$ is locally single-valued around $(y_0,x_0)$.
    
    Suppose first that $\Phi(x)=\varphi(x)+N_S(x)$. Let $z_0:=x_0+y_0-\varphi(x_0)$ and define
    $$
    N(z):=\varphi(\Pi_S(z))+(z-\Pi_S(z)),\quad z\in \bR^n.
    $$
    Consider the function $h:\gph\Phi\to\bR^n$ given by 
    $$
    h(x,y):=x+y-\varphi(x),\quad(x,y)\in\gph\Phi.
    $$
    {
    Since $y-\varphi(x)\in N_S(x)$ for any $(x,y)\in\gph\Phi$, we have $h(x,y)=z$ if and only if $x=\Pi_S(z)$ and $y=N(z)$. Thus, 
    $$
    h^{-1}(z)=(\Pi_S(z),N(z))\in\gph\Phi,\quad z\in\bR^n.
    $$
    The continuity of both $h$ and $h^{-1}$ shows that $h$ is a homeomorphism.
    
    We now show directly that $N^{-1}$ has the Aubin property around $(y_0,z_0)$. From the formula for $h^{-1}$ obtained above, we have
    $$
    N^{-1}(y)=\{h(x,y)\mid x\in\Phi^{-1}(y)\},\quad y\in\bR^n.
    $$
    Choose neighborhoods $U_0\subset\bR^n$ of $x_0$ and $V_0\subset\bR^n$ of $y_0$, and constants $\kappa,L>0$, such that, for all $y,y'\in V_0$,
    $$
    \Phi^{-1}(y)\cap U_0\subset \Phi^{-1}(y')+\kappa\|y-y'\|\mathbb B,
    $$
    and $h$ is $L$-Lipschitz on $(U_0\times V_0)\cap\gph\Phi$. Choose open neighborhoods $U_0^\prime\subset U_0$ of $x_0$ and $V_0^\prime\subset V_0$ of $y_0$ such that
    $$
    U_0^\prime+\kappa\|y-y'\|\mathbb B\subset U_0,\quad \forall y,y'\in V_0^\prime.
    $$
    By the continuity of $h^{-1}$ and $h^{-1}(z_0)=(x_0,y_0)$, choose an open neighborhood $W_0$ of $z_0$ such that $h^{-1}(W_0)\subset U_0^\prime\times V_0^\prime$. We claim that, for all $y,y'\in V_0^\prime$,
    $$
    N^{-1}(y)\cap W_0\subset N^{-1}(y')+L(\kappa+1)\|y-y'\|\mathbb B.
    $$
    Indeed, take any $z\in N^{-1}(y)\cap W_0$. Then $h^{-1}(z)=(x,y)$ for some $x\in\Phi^{-1}(y)\cap U_0^\prime$. By the Aubin property of $\Phi^{-1}$, there exists $x'\in\Phi^{-1}(y')$ such that $\|x-x'\|\le \kappa\|y-y'\|$. The choice of $U_0^\prime$ and $V_0^\prime$ gives $x'\in U_0$. Setting $z':=h(x',y')$, we have $z'\in N^{-1}(y')$, and
    $$
    \|z-z'\|=\|h(x,y)-h(x',y')\|\le L(\|x-x'\|+\|y-y'\|)\le L(\kappa+1)\|y-y'\|.
    $$
    }
    This proves the Aubin property of $N^{-1}$ around $(y_0,z_0)$. Applying Theorem \ref{NFGhom}, we obtain that $N$ is a local homeomorphism at $z_0$. It remains to translate this local invertibility back to $\Phi^{-1}$.
    
    There exists an open neighborhood $W\subset\bR^n$ of $z_0$ such that $N(W)$ is open and $N|_W:W\to N(W)$ is a homeomorphism. Since $h^{-1}(W)$ is a neighborhood of $(x_0,y_0)$ in $\gph\Phi$, there exist open neighborhoods $U\subset\bR^n$ of $x_0$ and $V\subset\bR^n$ of $y_0$, such that $(U\times V)\cap\gph\Phi\subset h^{-1}(W)$. Let
    $$
    W^\prime:=h((U\times V)\cap\gph\Phi),\qquad V^\prime:=N(W^\prime).
    $$
    {Since $h$ and $N|_W$ are homeomorphisms, $W'\subset W$ and $V'$ are open neighborhoods of $z_0$ and $y_0$, respectively.} Moreover, $V'\subset V$ follows from $h^{-1}(W')=(U\times V)\cap\operatorname{gph}\Phi$. Fix $y^\prime\in V^\prime$. Since $N|_W$ is a homeomorphism, there is a unique $z^\prime\in W^\prime$ such that $y^\prime=N(z^\prime)$. So $h^{-1}(z^\prime)=(\Pi_S(z^\prime),y^\prime)\in U\times V$. We claim that
    $$
    \Phi^{-1}(y^\prime)\cap U=\{\Pi_S(z^\prime)\}.
    $$
    Consider any $x^\prime\in\Phi^{-1}(y^\prime)\cap U$. Then $(x^\prime,y^\prime)\in(U\times V^\prime)\cap\gph\Phi$ and $h(x^\prime,y^\prime)\in W^\prime$. Since $N(h(x^\prime,y^\prime))=y^\prime=N(z^\prime)$, we have $h(x^\prime,y^\prime)=z^\prime$, and hence $x^\prime=\Pi_S(z^\prime)$. Therefore, the map $y\mapsto\Phi^{-1}(y)\cap U$ is single-valued on $V^\prime$, i.e., $\Phi^{-1}$ is locally single-valued around $(y_0,x_0)$.

    Now suppose that $\Phi(x)=\varphi(x)+N_S^{-1}(x)$. The argument above can be repeated with only minor modifications, which we record below. Set $N(z):=\varphi(z-\Pi_S(z))+\Pi_S(z)$ for $z\in\bR^n$ and $z_0:=x_0+y_0-\varphi(x_0)$, and define $h:\gph\Phi\to\bR^n$ by $h(x,y):=x+y-\varphi(x)$. Since $x\in N_S(y-\varphi(x))$ for any $(x,y)\in\gph\Phi$, we have $h(x,y)=z$ if and only if $x=z-\Pi_S(z)$ and $y=N(z)$. Thus $h^{-1}(z)=(z-\Pi_S(z),N(z))$ for all $z\in\bR^n$, so $h$ is a homeomorphism. Moreover, $N^{-1}(y)={h(x,y)\mid x\in\Phi^{-1}(y)}$ for all $y\in\bR^n$. Hence the preceding argument shows that $N^{-1}$ has the Aubin property around $(y_0,z_0)$, and Theorem \ref{NFGhom} implies that $N$ is a local homeomorphism at $z_0$. Repeating the final localization argument, with $x'=z'-\Pi_S(z')$ in place of $x'=\Pi_S(z')$, shows that $\Phi^{-1}$ is locally single-valued around $(y_0,x_0)$, completing the proof.
\end{proof}

{
\section{A counterexample without \texorpdfstring{$C^2$}{C2}-cone reducibility}
\label{sec:counter}

In this section, we construct a counterexample showing that the assumption of $C^2$-cone reducibility in Theorem \ref{main} cannot be dropped, even if $S$ is a three-dimensional semialgebraic set. This also explains, from another perspective, why Theorem \ref{NABhom} relies on a $C^2$-cone reduction to bring the problem into the framework of Theorem \ref{indexthm}, in which the reference point belongs to the set and the index is therefore shown to be $\pm 1$. In the absence of $C^2$-cone reducibility, this key reduction is unavailable, so the index is not necessarily $\pm 1$. Indeed, in our counterexample, the index of the associated map $N$ is $-2$, and strong regularity consequently fails.

Let $z=(z_1,z_2)\in\mathbb R^2$ and define $h:\mathbb R^2\to\mathbb R$ by
\begin{equation}
\label{eq:counter-h}
        h(z):=
        \begin{cases}
        \displaystyle
        \frac{z_1^3-3z_1z_2^2}{\sqrt{z_1^2+z_2^2}}, & {\rm if} \ z\ne0,\\
        0, & {\rm if} \ z=0.
        \end{cases}
\end{equation}
In fact, using polar coordinates $z=(r\cos\theta,r\sin\theta)$, we have
\begin{equation}
\label{eq:counter-h-polar}
        h(z)=r^2\cos 3\theta .
\end{equation}
This polar representation is used only to display the structure of $h$. In what follows, we work with the Cartesian coordinates $z=(z_1,z_2)$.
It is easy to check that the function $h$ defined in \eqref{eq:counter-h} is $C^1$, but it is not $C^2$ at $0$.
Set
\begin{equation}
\label{eq:counter-H}
        H(z):=\nabla h(z), \quad z\in\mathbb R^2.
\end{equation}
Then direct calculations give rise to   
\begin{equation}
\label{eq:counter-H-formula}
        H(0) = 0,
        \qquad
        H(z)
        =
        \frac{1}{\|z\|^3}
        \begin{pmatrix}
        2z_1^4+3z_1^2z_2^2-3z_2^4\\
        -z_1z_2(7z_1^2+3z_2^2)
        \end{pmatrix},\quad \forall z\ne 0.
\end{equation}
Define
\begin{equation}
\label{eq:counter-S}
        g(z):=4\|z\|^2+h(z),
        \qquad
        S:=\{(z,t)\in\mathbb R^2\times\mathbb R:\ t\ge g(z)\}.
\end{equation}
Finally, let
\begin{equation}
\label{eq:counter-Phi}
        f(z,t):=-4\|z\|^2,
        \qquad
        \Phi:=\nabla f+N_S,
        \qquad
        \bar x:=(0,0,0),
        \qquad
        \bar y:=(0,0,-1).
\end{equation}
Since $g(0)=0$ and $\nabla g(0)=0$, one has
\begin{equation*}
        N_S(\bar x)=\{(0,0,\lambda):\lambda\leq0\},
\end{equation*}
and therefore $\bar y\in\Phi(\bar x)$.

\begin{lemma}
The set $S$ defined in \eqref{eq:counter-S} is closed, convex and semialgebraic. Moreover, $S$ is not $C^2$-cone reducible at $\bar x$.
\end{lemma}

\begin{proof}
 Let $z$ be an arbitrary nonzero point in $R^2$ and define 
\begin{equation*}
\alpha(z):=\frac{z_1^3-3z_1z_2^2}{\|z\|^3}.
\end{equation*}
Then $|\alpha(z)|\le1$ and 
\begin{equation*}
        \operatorname{tr} H'(z)=-5\alpha(z),
        \qquad
        \det H'(z)=-9-5\alpha(z)^2.
\end{equation*}
Thus, the eigenvalues of $H'(z)=\nabla^2 h(z)$ are at least $-7$. Hence $\nabla^2 g(z)=8I+\nabla^2 h(z)\succeq I$. Since $g$ is $C^1$, it follows by integrating the above inequality along line segments that $g$ is strongly convex. Therefore $S=\operatorname{epi} g$ is closed and convex. The semialgebraicity of $S$ follows from the semialgebraicity of $g$.

It remains to show that $S$ is not $C^2$-cone reducible at $\bar x$. Suppose, to the contrary, that $S$ is $C^2$-cone reducible at $\bar x$. Then locally around $\bar x$,
\begin{equation*}
        S=\{x:\ \Xi(x)\in K\},
\end{equation*}
where $\Xi$ is $C^2$, $\Xi'(\bar x)$ is onto, and $K$ is a closed convex cone. We have $N_S(\bar x)=\nabla\Xi(\bar x) K^\circ$. Since $N_S(\bar x)$ is a ray and $\nabla\Xi(\bar x)$ is injective, $K^\circ$ must also be a ray. Hence $K$ is a halfspace. Consequently, near $\bar x$, the set $S$ is described by a single $C^2$ inequality with nonzero gradient. Thus $\partial S$ is a $C^2$ hypersurface near $\bar x$.

However, by \eqref{eq:counter-S}, $g$ is not $C^2$ at $0$ and $\partial S=\{(z,g(z)): z \in\bR^2\}$. This contradiction proves that $S$ is not $C^2$-cone reducible at $\bar x$.
\end{proof}

\begin{lemma}\label{ind2}
For the map $H$ defined in \eqref{eq:counter-H}, the inverse mapping $H^{-1}$ has the Aubin property around $(0,0)$. Moreover, $\operatorname{ind}(H,0)=-2$.
\end{lemma}

\begin{proof}
A direct computation from \eqref{eq:counter-H-formula} gives
\begin{equation}
\label{eq:counter-H-bound}
        2\|z\|\le \|H(z)\|\le3\|z\|,
        \quad \forall z\in\mathbb R^2.
\end{equation}
Moreover, $\det H'(z)<0$ for every $z\ne0$, $H(z)=H(-z)$, and $H^{-1}(0)=\{0\}$. In particular, $H$ is locally onto around $0$. By positive homogeneity, $H$ is onto $\mathbb R^2$.

We first compute the index. Let $e_1=(1,0)$ and $\mathbb{B}$ be the  closed unit ball in $\bR^2$. For $0<s<1$, the point $2se_1$ belongs to the same connected component of $\mathbb R^2\setminus H(\partial\mathbb B)$ as the origin, because $\|H(z)\|\ge2$ on $\partial\mathbb B$ by \eqref{eq:counter-H-bound}. Solving $H(z)=2se_1$ gives exactly two solutions $\pm se_1$.
Since $\det H'(z)<0$ at both solutions, degree theory gives
\begin{equation*}
        \operatorname{ind}(H,0)
        =
        \deg(H,\mathbb B,0)
        =
        \deg(H,\mathbb B,2se_1)
        =
        -2.
\end{equation*}

It remains to prove the Aubin property of $H^{-1}$, which is equivalent to the metric regularity of $H$. We claim that there exists $\kappa>0$ such that
\begin{equation*}
        \operatorname{dist}\bigl(z,H^{-1}(a)\bigr)
        \le
        \kappa\|H(z)-a\|,
        \quad
        \forall z,a\in\mathbb R^2.
\end{equation*}
Suppose not. Using the positive homogeneity of $H$, we may find sequences $\|z_k\|=1$ and $a_k\in\mathbb R^2$ such that
\begin{equation*}
        \operatorname{dist}\bigl(z_k,H^{-1}(a_k)\bigr)
        >
        k\|H(z_k)-a_k\|.
\end{equation*}
Since $H$ is onto, choose $\xi_k\in H^{-1}(a_k)$. By \eqref{eq:counter-H-bound},
\begin{equation*}
        \operatorname{dist}\bigl(z_k,H^{-1}(a_k)\bigr)
        \le \|z_k-\xi_k\|
        \le \|z_k\|+\|\xi_k\|
        \le 1+\frac{\|a_k\|}{2}.
\end{equation*}
Together with $\|H(z_k)\|\le3$, this implies that $\{a_k\}$ is bounded and
\begin{equation*}
        \|H(z_k)-a_k\|\to0.
\end{equation*}
Passing to a subsequence, assume that $z_k\to\hat z$ with $\|\hat z\|=1$. Then $a_k\to H(\hat z)$. Since $H'(\hat z)$ is nonsingular, $H$ has a Lipschitz local inverse near $\hat z$. Therefore, there exists $L>0$, such that for all large $k$,
\begin{equation*}
        \operatorname{dist}\bigl(z_k,H^{-1}(a_k)\bigr)
        \le
        L\|H(z_k)-a_k\|,
\end{equation*}
which is a contradiction. Thus $H^{-1}$ has the Aubin property around $(0,0)$.
\end{proof}

\begin{proposition}
For $S$, $f$, $\Phi$, $\bar x$ and $\bar y$ defined in \eqref{eq:counter-S} and \eqref{eq:counter-Phi}, $\Phi^{-1}$ has the Aubin property around $(\bar y,\bar x)$, but it is not strongly regular there. Moreover, let 
$$N(w):=\nabla f(\Pi_S(w))+(w-\Pi_S(w)),\quad w\in\bR^3.$$
We have $\ind(N,\bar y)=-2$.
\end{proposition}

\begin{proof}
Let $y=(p,\lambda)\in\mathbb R^2\times\mathbb R$ such that $\lambda<0$, and $x=(z,t)\in\mathbb R^2\times\mathbb R$ such that $y\in\Phi(x)$. Then $x$ cannot lie in the interior of $S$, otherwise the third component of $\Phi(x)$ equals $0$. Hence $x\in\partial S$. By \eqref{eq:counter-S}, $y=(-8z,0)+(-\lambda\nabla g(z),\lambda)$.
Since $\nabla g(z)=8z+H(z)$, this is equivalent to
\begin{equation}
\label{eq:counter-T}
        y=T(z,\lambda)
        :=
        (H(z)-(\lambda+1)(8z+H(z)),\lambda).
\end{equation}
Set
\begin{equation*}
        T_0(z,\lambda):=(H(z),\lambda).
\end{equation*}
By Lemma \ref{ind2}, $H^{-1}$ has the Aubin property around $(0,0)$, and hence $T_0^{-1}$ has the Aubin property around $(\bar y,\bar y)$. Moreover,
\begin{equation*}
        T(z,\lambda)-T_0(z,\lambda)
        =
        ((\lambda+1)(8z+H(z)),0).
\end{equation*}
Lemma \ref{ss} ensures that this perturbation is strictly stationary at $\bar y$. So, Lemma \ref{invmapthm} implies that $T^{-1}$ has the Aubin property around $(\bar y,\bar y)$. Finally, the localized inverse of $\Phi$ around $(\bar y,\bar x)$ is obtained from $T^{-1}$ through the Lipschitz map
\begin{equation*}
        (z,\lambda)\mapsto (z,g(z)).
\end{equation*}
Consequently, $\Phi^{-1}$ has the Aubin property around $(\bar y,\bar x)$.

By Lemmas \ref{keepindex} and \ref{ind2}, we have $\ind(T,\bar y)=\ind(T_0,\bar y)=\ind(H,0)=-2$. Let
$$Q(z,\lambda):=(z,g(z))+\lambda(-\nabla g(z),1),\quad\lambda<0.$$
Then $\Pi_S(Q(z,\lambda))=(z,g(z))$. Thus, $N\circ Q=T$ around $\bar y$. Moreover, $Q$ is a local homeomorphism at $\bar y$ with $\ind(Q,\bar y)=1$. Therefore, \cite[Theorem 7.4.3]{dinca2021brouwer} yields $\ind(N,\bar y)=\ind(T,\bar y)=-2$.

It remains to show that strong regularity fails. For $s>0$, define $y_s:=(2s,0,-1)$.
Let $e_1=(1,0)$. By \eqref{eq:counter-H-formula} and \eqref{eq:counter-T}, we obtain
\begin{equation*}
        T(se_1,-1)=T(-se_1,-1)=y_s.
\end{equation*}
Equivalently, $(se_1,g(se_1))$ and $(-se_1,g(-se_1))$ are two distinct points of $\Phi^{-1}(y_s)$. Both converge to $\bar x$ as $s\to0$, while $y_s\to\bar y$. Therefore, $\Phi^{-1}$ is not locally single-valued around $(\bar y,\bar x)$, thus is not strongly regular around $(\bar y,\bar x)$.
\end{proof}

\begin{remark}
    In fact, $\Phi^{-1}$ is locally two-valued around $(\bar y,\bar x)$. The number two is not essential. For any integer $n\geq1$, replacing $h$ by $h_n(z)=r^2\cos((n+1)\theta)$ and replacing the coefficient $4$ in $g$ and $f$ by a sufficiently large constant $c_n$, one obtains the same type of counterexample. Indeed, the index of $\nabla h_n$ at $0$ equals to $-n$, and the localized inverse solution mapping has $n$ distinct branches.
\end{remark}

\section{Application to KKT systems}
\label{sec:app}

Consider the constrained optimization problem
\begin{equation}
\label{eq:app-op}
    \min_{x\in\bR^n} f(x)\quad\text{s.t.}\quad G(x)\in \mathcal K,
\end{equation}
where $f:\bR^n\to\bR$ and $G:\bR^n\to\bR^m$ are twice continuously
differentiable functions and $\mathcal K\subset\bR^m$ is a closed convex set. The Lagrangian function of \eqref{eq:app-op} is defined by
$$\mathcal L(x,y):=f(x)+\langle y,G(x)\rangle,\quad(x,y)\in\bR^n\times\bR^m.$$
And the KKT system of \eqref{eq:app-op} is
\begin{equation}
\label{eq:app-kkt}
    0=\nabla_x\mathcal L(x,y)\quad\text{and}\quad y\in N_{\mathcal K}(G(x)).
\end{equation}
For any solution $(\bar x,\bar y)$ to the KKT system \eqref{eq:app-kkt}, $\bar x$ is called a stationary point of \eqref{eq:app-op}, and $\bar y$ is called a (Lagrange) multiplier associated with $\bar x$.

The canonical perturbation \cite[Section 5.1.3]{bonnans2000perturbation} of \eqref{eq:app-op} is given by
\begin{equation}
\label{eq:app-op-p}
    \min_{x\in\bR^n} f(x)-\langle a,x\rangle\quad\text{s.t.}\quad G(x)+b\in \mathcal K,
\end{equation}
where $a\in\bR^n$ and $b\in\bR^m$ are perturbation vectors. The perturbed version of \eqref{eq:app-kkt}, i.e., the KKT system of \eqref{eq:app-op-p} is given by
\begin{equation}
\label{eq:app-kkt-p}
    a=\nabla_x\mathcal L(x,y)\quad\text{and}\quad y\in N_{\mathcal K}(G(x)+b).
\end{equation}
We reformulate the perturbed KKT system \eqref{eq:app-kkt-p} as a generalized equation of the form in \eqref{ge}:
\begin{equation}
\label{eq:app-kkt-phi}
    \begin{pmatrix}
        a \\ b
    \end{pmatrix}
    \in
    \Phi_\text{KKT}(x,y):=
    \begin{pmatrix}
        \nabla_{x}{\mathcal L}(x,y) \\ -G(x)
    \end{pmatrix}+
    N_{\{0\}\times\mathcal K}^{-1}
    \begin{pmatrix}
        x \\ y
    \end{pmatrix}.
\end{equation}

Let $(\bar x,\bar y)\in \Phi_\mathrm{KKT}^{-1}(0,0)$ be a solution of \eqref{eq:app-kkt}. A central issue is to understand the local behavior of the solution mapping $\Phi_\mathrm{KKT}^{-1}$ around $((0,0),(\bar x,\bar y))$, particularly its Aubin property and strong regularity. A fundamental result due to Dontchev and Rockafellar \cite{dontchev1996characterizations} states that, when $\mathcal K$ is a polyhedral set, the Aubin property and strong regularity for $\Phi_\mathrm{KKT}^{-1}$ are equivalent. This equivalence has also been established when $\mathcal K$ is a Cartesian product of second-order cones \cite{chen2025aubin} or a Cartesian product of positive semidefinite cones \cite{chen2025characterizations}, under the assumption that $\bar x$ is a locally optimal solution of \eqref{eq:app-op}. The presence of this additional local optimality assumption distinguishes the latter two results from the classical polyhedral one. With Theorem \ref{main}, we can remove the local optimality assumption and derive  the equivalence for a much broader class of sets $\mathcal K$.

\begin{theorem}
\label{thm:app-kkt}
    Let $(\bar x,\bar y)$ be a solution to the KKT system of \eqref{eq:app-op}. Suppose that the set $\mathcal K$ in \eqref{eq:app-op} is either a $C^2$-cone reducible set or a cone whose polar is $C^2$-cone reducible. Then the Aubin property and strong regularity for the solution mapping $\Phi_\mathrm{KKT}^{-1}$ in \eqref{eq:app-kkt-phi} are equivalent around $((0,0),(\bar x,\bar y))$.
\end{theorem}

\begin{proof}
    If $\mathcal K$ is a $C^2$-cone reducible set, then the desired conclusion follows directly from the second form of Theorem \ref{main}. It remains to consider the case where $\mathcal K$ is a cone and $\mathcal K^\circ$ is $C^2$-cone reducible. In this case, using the identity
    $$N^{-1}_{\{0\}\times\mathcal K}=N_{\bR^n\times\mathcal K^\circ},$$
    the equivalence follows from the first form of Theorem \ref{main}.
\end{proof}

The class of sets $\mathcal K$ covered by Theorem \ref{thm:app-kkt} is substantially broader, encompassing all the cases mentioned above and including, in particular, the $p$-order cone, which is the epigraph of the $p$-norm
\begin{equation*}
    \mathcal K_p:=\operatorname{epi}\|\cdot\|_p
    =\{(t,z)\in\mathbb R\times\mathbb R^r\mid t\ge\|z\|_p\}.
\end{equation*}
When $p=2$, it becomes second-order cone.  

\begin{corollary}
The equivalence in Theorem \ref{thm:app-kkt} holds, in particular, when $\mathcal{K}$ is a Cartesian product of polyhedral cones, positive semidefinite cones, and $p_i$-order cones, with either all $p_i\in[2,\infty)$ or all $p_i\in(1,2]$.
\end{corollary}

\begin{proof}
For $p\in[2,\infty)$, the $p$-order cone is $C^2$-cone reducible: away from the origin, it is locally represented by a $C^2$ inequality with nonzero gradient on the boundary point, while at the origin the reduction is trivial. Hence, if all $p_i\in[2,\infty)$, then every factor in $\mathcal K$ is $C^2$-cone reducible. Thus, $\mathcal K$ is also $C^2$-cone reducible, and the conclusion follows from Theorem~\ref{thm:app-kkt}.

If all $p_i\in(1,2]$, let $q_i=p_i/(p_i-1)\in[2,\infty)$. Then
$$
\mathcal K_{p_i}^{\circ}=-\mathcal K_{q_i},
$$
so the polar of each $p_i$-order cone is $C^2$-cone reducible. Moreover, polars of polyhedral cones and polars of positive semidefinite cones are $C^2$-cone reducible. Therefore, $\mathcal K^\circ$ is $C^2$-cone reducible, and the conclusion follows from Theorem \ref{thm:app-kkt}.
\end{proof}
}

\section{Conclusions}\label{sec:con}

In this paper, we establish the equivalence between the Aubin property and strong regularity for generalized equations over $C^2$-cone reducible sets. This result resolves a long-standing open question in variational analysis and extends the classical theorem of Dontchev and Rockafellar \cite{dontchev1996characterizations} beyond the polyhedral setting to a broad class of non-polyhedral sets.

{
Our proof  strategy departs from traditional variational techniques by combining convex geometry with the theory of topological degrees.  We also constructed a three-dimensional closed convex semialgebraic set that is not $C^2$-cone reducible, together with a generalized equation whose solution mapping has the Aubin property but is not strongly regular. This shows that the $C^2$-cone reducibility assumption cannot be removed in general.

Finally, we applied the main result to KKT systems for general constrained optimization problems. Without assuming local optimality, we proved that the Aubin property and the strong regularity of the KKT solution mapping are equivalent at any KKT point when the constraint set is $C^2$-cone reducible or is a cone whose polar is $C^2$-cone reducible. This provides a unified framework for analyzing the stability of important optimization problems, including conventional nonlinear programming, nonlinear second-order cone programming, and nonlinear semidefinite programming. Moreover, the result covers broader classes of conic constraints, including suitable Cartesian products of polyhedral cones, positive semidefinite cones, and $p$-order cones.
}

\bibliographystyle{apalike}
\bibliography{ref}

\end{document}